\documentclass[11pt]{article}
\usepackage{amsfonts}
\def\-{\hbox{-}}
\def\.{{\,\hskip-1pt\cdot\hskip-1pt\,}}
\def\O{{\cal O}}
\def\i{\leq}
\def\K{{\cal K}}

\def\E{{\cal E}}

\def\T{{\cal T}}

\def\A{{\cal A}}
\def\B{{\cal B}}

\def\Z{{\cal Z}}

\def\Ab{\frak A\frak b}

\def\int{\frak i\frak n\frak t}

\def\qq{\quad{\rm and}\quad}

\def\too{\longrightarrow}

\input amssym.def
\input amssym.tex
\usepackage[all]{xy}
\usepackage{eufrak}
\usepackage{amscd}

\begin{document}
\begin{center}{\large\bf
Glauberman correspondents and extensions of nilpotent block
algebras}

\bigskip\large\bf Lluis Puig and Yuanyang Zhou

\end{center}

\insert\footins{\scriptsize 2000 {\it Mathematics Subject
Classification}. Primary 20C15, 20C20

\smallskip The second author is supported by Program for New Century
Excellent Talents in University and by NSFC (No. 11071091).}

\begin{abstract} The main purpose of this paper is to prove that the extensions of a nilpotent
block algebra and its Glauberman correspondent block algebra are
Morita equivalent under an additional group-theoretic condition (see
Theorem 1.6); in particular, Harris and Linckelman's theorem and Koshitani and
Michler's theorem are covered (see Theorems £7.5 and £7.6). The ingredient to carry out our  purpose
is the two main results in K\"ulshammer and Puig's work {\it Extensions of nilpotent blocks\/}; we actually revisited
them, giving completely new proofs of both and slightly improving the second one
(see Theorems £3.5 and £3.14).

\end{abstract}

\bigskip\noindent{\bf\large 1. Introduction}

\bigskip\noindent{\bf 1.1.}\quad
Let $\O$ be a complete discrete valuation ring with an algebraically
closed residue field $k$ of characteristic $p$ and a quotient field
$\K$ of characteristic 0. In addition, $\K$ is also assumed to be
big enough for all finite groups that we consider below. Let $H$ be
a finite group. We denote by ${\rm Irr}_\K(H)$ the set of all
irreducible characters of $H$ over $\K$.
Let $A$ be another finite group and assume that there is a group homomorphism $A\rightarrow {\rm Aut}(H)$. Such a group $H$ with an $A$-action is called an $A$-group. We denote by $H^A$ the subgroup of all
$A$-fixed elements in $H$. Clearly $A$ acts on ${\rm Irr}_\K(H)$. We
denote by ${\rm Irr}_{\cal K}(H)^A$ the set of all $A$-fixed
elements in ${\rm Irr}_{\cal K}(H)$. Assume that $A$ is solvable and
the order of $A$ is coprime to the order of $H$. By
\cite[Theorem 13.1]{I}, there is a bijection
$$\pi(H,\,A): {\rm Irr}_{\cal K}(H)^A \rightarrow
{\rm Irr}_{\cal K}(H^A)$$ such that

\medskip\noindent{\bf 1.1.1.} For any normal subgroup $B$ of $A$, the
bijection $\pi(H,\, B)$ maps ${\rm Irr}_{\K}(H)^A$ to
${\rm Irr}_{\K}(H^B)^A$, and in ${\rm Irr}_{\K}(H)^A$ we have
$$\pi(H,\,A) = \pi(H^B,\,A/B) \circ \pi(H,\,B)\,.$$

\medskip\noindent{\bf 1.1.2.}  If $A$ is a $q$-group for some prime $q$, then for
any $\chi\in {\rm Irr}_{\K}(H)^A$, the corresponding irreducible
character $\pi(H,\,A)(\chi)$ of $G^A$ is the unique irreducible
constituent of ${\rm Res}^H_{H^A}(\chi)$ occurring with a
multiplicity coprime to~$q$.

\medskip\noindent The character $\pi(H,\,A)(\chi)$ of $H^A$ is called the {\it Glauberman correspondent\/} of the character $\chi$ of $H$.

\bigskip\noindent{\bf 1.2.}\quad
For any central idempotent $c$ of $\O H$, we denote by ${\rm
Irr}_{\K}(H, c)$ the set of all irreducible characters of $H$
provided by some $\K Hc$-module. Let $b$ be a block of $H$ --- namely $b$ is a primitive central idempotent of $\O H\,;$ then $\O Hb$ is called
the {\it  block algebra\/} corresponding to $b$. Assume that $A$ stabilizes
the block $b$ and  centralizes a defect group of $b$.
Then, by \cite[Proposition 1 and Theorem 1]{W}, $A$ stabilizes all characters
of ${\rm Irr}_\K(H, b)$ and there is a unique block ${\it w}(b)$ of
$\O (H^A)$ such that
$${\rm Irr}_{\K}(H^A, {\it w}(b))=\pi(H,\,A)({\rm
Irr}_\K(H, b))\,;$$
moreover ,there is a perfect isometry (see \cite{B1})
$$R_H^b: {\cal R}_\K (H, b)\rightarrow {\cal R}_\K (H^A, {\it w}(b))$$
such that $R_H^b(\chi)=\pm\pi(H,\,A)(\chi)$ for any $\chi\in {\rm
Irr}_{\K}(H, b)$, where we denote by ${\cal R}_\K (H, b)$ and ${\cal
R}_\K (H^A, {\it w}(b))$ the additive groups generated by ${\rm
Irr}_{\K}(H, b)$ and ${\rm Irr}_{\K}(H^A, {\it w}(b))$. Such a block
${\it w}(b)$ is called the {\it Glauberman correspondent\/} of $b$ (see
\cite{W}). Since a perfect isometry between blocks is  often nothing but
 the character-theoretic `shadow' of a derived equivalence, it seems reasonable
 to ask whether there is a derived equivalence between a block and its Glauberman correspondent. In the last few years, some Morita equivalences between $b$ and ${\it w}(b)$ were found in the cases where $H$ is $p$-solvable or the defect
 groups of $b$ are normal in $H$, which supply Glauberman correspondences
 from ${\rm Irr}_{\K}(H, b)$ to ${\rm Irr}_{\K}(H^A, {\it w}(b))$ (see \cite{H},
\cite{KG} and \cite{H1}); moreover, all these Morita equivalences between $b$
and ${\it w}(b)$ are {\it basic} in the sense of~\cite{P1}.

\bigskip\noindent{\bf 1.3.}\quad By induction, the groups $H$ and
$H^A$ and the blocks $b$ and ${\it w}(b)$ in the main results of
\cite{H}, \cite{KG} and \cite{H1} can be reduced to the situation where,
for some $A$-stable normal subgroup $K$ of $H\,,$ we have
$H=H^A \.   K$
, the block $b$ is an
$H$-stable block of $K$ with trivial or central defect group, and
the block ${\it w}(b)$ is an $H^A$-stable block of $K^A$ with trivial
or central defect group. Recall that the block $b$ of $H$ is called {\it nilpotent\/}
(see \cite{P4}) if the quotient group $N_H(R_\varepsilon)/C_H(R)$ is a $p$-group for any local pointed group $R_\varepsilon$ on $\O H b$.
Blocks with trivial or central defect
group are nilpotent and therefore in these situations $\O Hb$ and $\O (H^A){\it w}(b)$ are extensions of the nilpotent block algebras $\O K b$
and $\O K^A{\it w}(b)$ respectively. K\"ulshammer and Puig already
precisely described the algebraic structure of extensions of
nilpotent block algebras (see \cite{KP} or Section 3 below) and these results can be applied to blocks of $p$-solvable groups (see \cite{P2}) and to blocks with normal
defect groups (see \cite{R, K}). Thus, it is reasonable to seek a common generalization of the main results of
\cite{H, KG, H1} in the setting of extensions of nilpotent block
algebras.

\bigskip\noindent{\bf 1.4.}\quad Let $G$ be another finite $A$-group having $H$ as
an $A$-stable normal subgroup and consider the $A$-action on $H$
induced by the $A$-group $G$. We assume that $A$ stabilizes $b$ and denote by $N$ the stabilizer of $b$ in $G$.
Clearly $N$ is $A$-stable. Set
\begin{center}$
c={\rm Tr}_N^G (b)$ and $\alpha=\{c\}$\quad ;
\end{center}
then the idempotent $c$ is $A$-stable and $\alpha$ is an
$A$-stable point of $G$ on the group algebra $\O H$ (the action of
$G$ on $\O H$ is induced by conjugation). In particular, $G_\alpha$
is a pointed group on $\O H$. Let~$P$ be a defect group of
$G_\alpha\,;$ then, by \cite[Proposition 5.3]{KP}, $Q=P\cap H$ is a defect
group of the block $b$ of~$H$.

\bigskip\noindent{\bf Theorem 1.5.}\quad {\it Assume that $A$ centralizes $P$,
that $A$ is solvable and that the orders of $G$ and $A$ are coprime. Set ${\it
w}(c)={\rm Tr}_{N^A}^{G^A}({\it w}(b))$ and ${\it
w}(\alpha)=\{{\it w}(c)\}$. Then, ${\it w}(\alpha)$ is a point of $G^A$ on the group algebra $\O (H^A)$ and $P$ is a defect group of the pointed group
$(G^A)_{{\it w}(\alpha)}$ on $\O (H^A)\,.$ Moreover, if $G=H\.    G^A$ and the block $b$ of $H$ is nilpotent, we have
\begin{center}${\rm Irr}_{\K}(G, c)={\rm Irr}_{\K}(G, c)^A$ and $\pi(G, A)({\rm Irr}_{\K}(G, c))={\rm Irr}_{\K}(G^A, {\it w}(c))$. \end{center}}

\bigskip The following theorem shows that there is a ``{\it basic\/}" Morita
equivalence between $\O G c$ and $\O G^A {\it w}(c)$; that is to say, this Morita equivalence induces basic Morita equivalences~\cite{P1} between corresponding block algebras.

\bigskip\noindent{\bf Theorem 1.6.}\quad {\it Assume that $A$ centralizes $P$,
that $A$ is solvable and that the orders of $G$ and $A$ are coprime. Set ${\it
w}(c)={\rm Tr}_{N^A}^{G^A}({\it w}(b))$. Assume that $G=G^A\.     H$
and that the block $b$ is nilpotent. Then, there is an $\O (H\times
H^A)$-module $M$ inducing a basic Morita equivalence
between $\O Hb$ and $\O (H^A) {\it w}(b)\,,$ which can be extended to the inverse image $K$ in $N\times N^A$ of the ``diagonal'' subgroup of $N/H\times N^A/H^A$
in such a way that ${\rm Ind}^{G\times G^A}_K (M)$
induces a ``basic" Morita equivalence between $\O Gc$ and $\O (G^A){\it w}(c)$.
}

\bigskip\noindent{\bf Remark 1.7.}\quad Since $G=H\.     G^A$, we have $N=H\.     N^A$ and
then the inclusion $N^A\subset N$ induces a group isomorphism
$N/H\cong N^A/H^A$.

\bigskip We use pointed groups introduced by Lluis Puig. For more details on pointed groups,
readers can see \cite{P5} or Paragraph 2.5 below. In Section 2, we
introduce some notation and terminology. Section 3 revisits
K\"ulshammer and Puig's main results on extensions of nilpotent
blocks; the proof of the existence and uniqueness of the finite group $L$ (see
 \cite[Theorem 1.8]{KP} and Theorem 3.5 below) is  dramatically simplified; actually, Corollary 3.14 below slightly improves  \cite[Theorem 1.12]{KP}; explicitly,
$S_\gamma$ in Corollary 3.14 is unique up to determinant one.
With the Glauberman correspondents of
blocks due to Watanabe, in Section 4 we define Glauberman
correspondents of extensions of blocks and compare the local
structures of extensions of blocks and their Glauberman
correspondents.

\medskip
By Puig's structure theorem of nilpotent blocks,
there is a bijection between the sets of irreducible characters of the nilpotent
block $b$ of $H$ and of its defect $Q$; in Section 5, for a suitable local point
$\delta$ of $Q\,,$ we prove that
this bijection preserves $N_G(Q_\delta)$-actions on these sets. As a consequence,
we obtain an $N_G(Q_\delta)$-stable
irreducible character $\chi$ of $H$ such that $\chi$ lifts the
unique irreducible Brauer character of the nilpotent block $b$ of
$H$ and that the Glauberman correspondent character $\pi(H, A)(\chi)$ lifts the
unique irreducible Brauer character of the Glauberman correspondent block
${\it w}(b)$ of $H^A$ (see Lemma 5.6).

\medskip
Obviously, $N$ stabilizes the
unique simple module in the nilpotent block $b$ of $H$; with this
$N$-stable $\O H b$-simple module, we construct an $A$-stable
$k^*$-group $\skew3\hat {\bar N}^{^k}$ (see £2.3 and £3.13 below);
since $N^A$ stabilizes the unique simple module of the nilpotent block
${\it w}(b)$ of $H^A\,,$  a $k^*$-group $\,\widehat{\overline{\! N^A}}^k$ is similarly
constructed. In Section 6, we prove that $\,\widehat{\overline{\! N^A}}^k$
and $(\skew3\hat{\bar N}^{^k})^A$ are isomorphic as $k^*$-groups (see Theorem
6.4). In Section 7, we use the improved version of K\"ulshammer and
Puig's main result to prove our main theorem 1.6.
\vfill
\eject

\vskip 1cm \noindent{\bf\large 2. Notation and terminology}

\bigskip\noindent{\bf 2.1.}\quad Throughout this paper, all
$\O$-modules are  $\O$-free finitely generated --- except in 2.4 below; all $\O$-algebras
have identity elements, but their subalgebras need not have the
same identity element. Let $\cal A$ be an
$\O$-algebra; we denote by $\A^\circ\,,$ ${\cal A}^*$, $Z({\cal A})$, $J({\cal
A})$ and $1_{\cal A}$ the opposite $\O\-$algebra of $\A\,,$ the multiplicative
group of all invertible elements of ${\cal A}$, the center of~${\cal A}$, the radical of
${\cal A}$ and the identity element of ${\cal A}$ respectively.
Sometimes we write $1$ instead of $1_{\cal A}\,.$ For any abelian group $V$, ${\rm id }_V$ denotes the identity automorphism on $V$. Let
${\cal B}$ be an $\O$-algebra; a homomorphism ${\cal F}: {\cal
A}\rightarrow {\cal B}$ of $\O$-algebras is said to be an {\it embedding\/}
if ${\cal F}$ is injective and we have
 $${\cal F}({\cal A})={\cal F}(1_{\cal A}){\cal B}{\cal F}(1_{\cal
A})\quad .$$
 Let $S$ be a set and $G$ be a group acting on $S$. For any
$g\in G$ and $s\in S$, we write the action of $g$ on $s$ as $s\.g$.

\bigskip\noindent{\bf 2.2.}\quad Let $X$ be a finite group. An
$X$-interior $\O$-algebra ${\cal A}$ is an $\O$-algebra ${\cal A}$ endowed
with a group homomorphism $\rho:X\rightarrow {\cal A}^*$; for any
$x, y\in X $ and $a\in {\cal A}$, we write $\rho(x)a\rho(y)$ as $x\.     a\.     y$ or $xay$ if there is no confusion. Let $\varrho: Y\rightarrow X$ be a group
homomorphism; the $\O$-algebra ${\cal A}$ with the group
homomorphism $\rho\circ\varrho: Y\rightarrow {\cal A}^*$ is an $Y$-interior $\O$-algebra and we denote it by ${\rm Res}_{\varrho}({\cal
A})$. Let ${\cal A}'$ be another $X$-interior $\O$-algebra; an
$\O$-algebra homomorphism ${\cal F}:{\cal A}\rightarrow {\cal A}'$ is said
to be a homomorphism of $X$-interior $\O$-algebras if for any $x, y\in X $ and
any $a\in {\cal A}$, we have
${\cal F}(xay)=x{\cal F}(a)y$. The tensor product ${\cal A}\bigotimes_\O {\cal
A}'$ of ${\cal A}$ and ${\cal A}'$ is an $X$-interior $\O$-algebra
with the group homomorphism
$$ X\rightarrow ({\cal A}\otimes_\O
{\cal A}')^*,\quad x\mapsto x1_{\cal A}\otimes x1_{{\cal A}'}\quad
.$$ Let $Z$ be a subgroup of $X$ and let ${\cal B}$ be an $\O
Z$-interior algebra. Obviously, the left and right multiplications by
$\O Z$ on ${\cal B}$ define an $(\O Z, \O Z)$-bimodule structure on
${\cal B}$. Set
$${\rm Ind}_Z^X ({\cal B})=\O X\otimes_{\O Z} {\cal B}\otimes_{\O Z} \O X
$$
and then this the $(\O X, \O X)$-bimodule ${\rm Ind}_Z^X ({\cal B})$ becomes an
$X$-interior $\O$-algebra with the product
$$(x\otimes b\otimes y)(x'\otimes b'\otimes y')
=\cases{x\otimes b\. y x'\.b'\otimes y'&if $yx'\in Z$\cr {}&{}\cr 0
&otherwise\cr}$$ for any $x, y, x', y'\in X$ and any $b,b'\in {\cal
B}\,,$ and with the homomorphism $\O X\rightarrow {\rm Ind}_Z^X
({\cal B})$ mapping $x\in X$ onto $\sum_y xy\otimes 1\otimes
y^{-1}$, where $y$ runs over a set of representatives for left cosets
of $Z$ in $X$.

\bigskip\noindent{\bf 2.3.}\quad A $k^*$-group
with $k^*$-quotient $X$ is a group $\hat X$ endowed with an injective group
homomorphism $\theta: k^*\rightarrow Z(\hat X)$ together with an isomorphism
$\hat X/\theta(k^*)\cong X$; usually we omit to mention $\theta$ and the quotient $X=\hat X/\theta(k^*)$ is called the $k^*$-quotient of $\hat X$, writing $\lambda\.     \hat x$ instead of $\theta(\lambda)\hat x$ for any $\lambda\in k^*$ and any $\hat x\in \hat X$.
We denote by $\hat Y$ the inverse image of $Y$ in $\hat X$ for
any subset $Y$ of~$X$ and, if  no precision is needed, we often denote by $\hat x$ some lifting of an element $x\in X$.
We denote by $\hat X^\circ$ the $k^*$-group with the same underlying group  $\hat X$ endowed with the group homomorphism
$\theta^{-1}: k^*\rightarrow Z(\hat X),\, \lambda\mapsto \theta(\lambda)^{-1}$. Let $\vartheta: Z\rightarrow X$ be a group homomorphism; we denote by ${\rm Res}_{\vartheta}(\hat X)$ the $k^*$-group formed by the group of pairs
$(\hat x, y)\in \hat X\times Z$ such that $\vartheta(y)$ is the
image of $\hat x$ in $X$, endowed with the group homomorphism mapping $\lambda\in k^*$ on $(\theta(\lambda), 1)$; up to suitable identifications, $Z$ is the $k^*$-quotient of ${\rm
Res}_{\vartheta}(\hat X)$. Let $\hat U$ be another $k^*$-group with $k^*$-quotient
$U$. A group homomorphism $\phi: \hat X\rightarrow \hat U$ is a
homomorphism of $k^*$-groups if $\phi(\lambda\.     \hat x)=\lambda\.     \phi(\hat x)$
for any $\lambda\in k^*$ and $\hat x\in \hat X$. For more details on
$k^*$-groups, please see \cite[\S 5]{P6}.

\bigskip\noindent{\bf 2.4.}\quad Let $\hat X$ be a $k^*$-group with $k^*$-quotient $X$. By \cite[Charpter II, Proposition 8]{S}, there exists a
canonical decomposition $\O^*\cong (1+J(\O))\times k^*$, thus $k^*$
can be canonically regarded as a subgroup of $\O^*$. Set
$$\O_*\hat
X=\O\otimes_{\O k^*}\O \hat X
\quad, $$
where the left $\O k^*$-module $\O\hat X$ and the right $\O k^*$-module
$\O$ are defined by the left and right multiplication by $k^*$ on $\hat X$ and $\O
^*$ respectively. It is straightforward to verify that the
$\O$-module  $\O_*\hat X$ is an $\O$-algebra with the distributive
multiplication
$$(a_1\otimes\hat x_1)(a_2\otimes\hat {x}_2)
=a_1a_2\otimes\hat x_1\hat{x}_2$$ for any $a_1,a_2\in \O$ and any $\hat
x_1,\hat{x}_2 \in \hat X$.

\bigskip\noindent{\bf 2.5.}\quad Let ${\cal A}$ be an
$X$-algebra over $\O$; that is to say, $\A$ is endowed with a group homomorphism
$\psi: X\rightarrow {\rm Aut}({\cal A})$, where ${\rm Aut}({\cal A})$ is
the group of all $\O$-automorphisms of $A\,;$ usually, we omit to mention $\psi\,.$ For any subgroup $Y$ of
$X$, we denote by ${\cal A}^Y$ the $\O$-subalgebra of all $Y$-fixed
elements in ${\cal A}$. A {\it pointed group\/} $Y_\beta$ on ${\cal A}$
consists of a subgroup $Y$ of $X$ and of an $({\cal A}^Y)^*$-conjugate
class $\beta$ of primitive idempotents of ${\cal A}^Y$. We often say that
$\beta$ is a {\it point\/} of $Y$ on ${\cal A}$. Obviously, $X$ acts on the
set of all pointed groups on ${\cal A}$ by the equality $(Y_\beta)^x=Y^x_{\psi(x^{-1})(\beta)}$ and we denote by
$N_X(Y_\beta)$ the stabilizer of $Y_\beta$ in $X$ for any pointed
group $Y_\beta$ on ${\cal A}$. Another pointed group $Z_\gamma$ is said
{\it contained in\/} $Y_\beta$ if $Z\leq Y$ and there exist some
$i\in \beta$ and $j\in \gamma$ such that $ij=ji=j$. For a subgroup $U$ of $G$, set
$${\cal A}(U)=k\otimes _\O ({\cal A}^U/\sum_V {\cal A}^U_V)\quad$$ where $V$ runs over the set of proper subgroups of $U$ and ${\cal A}^U_V$ is the image
of the relative trace map ${\rm Tr}_V^U: {\cal A}^V\rightarrow {\cal
A}^U$; the canonical surjective homomorphism ${\rm Br}^{\cal A}_U:
{\cal A}^U\rightarrow {\cal A}(U)$ is called the {\it Brauer homomorphism\/} of the
$X$-algebra ${\cal A}$ at $U$. When ${\cal A}$ is equal to the group
algebra $\O X$, the homomorphism $kC_X(U)\rightarrow {\cal A}(U)$
sending $x\in C_X(U)$ onto the image of $x$ in ${\cal A}(U)$ is an
isomorphism, through which we identify ${\cal A}(U)$ with
$kC_X(U)$. A pointed group $U_\gamma$ on $\cal A$ is said {\it local\/} if the
image of $\gamma$ in ${\cal A}(U)$ is not equal to $\{0\}\,,$ which
forces $U$  to be a $p$-group; then, a local pointed group $U_\gamma$
is said a {\it defect pointed group\/} of a pointed group $Y_\beta$ on $\cal A$
if $U_\gamma\leq Y_\beta$ and we have $\beta\subset {\rm Tr}_U^Z({\cal
A}^U\.    \gamma\.     A^U)$, where ${\cal A}^U\.    \gamma\.     A^U$
is the ideal of ${\cal A}^U$ generated by~$\gamma$. Let $c$ be a
block of $X\,;$ then $\{c\}$ is a point of $X$ on $\O X$ and if
$P_\gamma$ is a defect pointed group of $X_{\{c\}}$ then $P$ is a
defect group of $c$.

\vskip 1cm \noindent{\bf\large 3. Extensions of nilpotent blocks
revisited}

\bigskip In this section, we assume that $\O$ is a complete
discrete valuation ring with an algebraically closed residue
field of characteristic $p$.

\bigskip
\noindent {\bf £3.1.}\quad Let $G$ be a finite group, $H$ be a
normal subgroup of $G$ and $b$ be a block of $H$ over~$\O$. Denote by $N$
the stabilizer of $b$ in $G$ and set $\bar N=N/H$. Obviously,
$\beta=\{b\}$ is a point of $H$ and $N$ on $\O H$ and there is a unique
pointed group $G_\alpha$ on $\O H$ such that
$$H_\beta\leq N_\beta\leq G_\alpha\quad .$$
Let $Q_\delta$ be a defect pointed group of~$H_\beta$ and $P_\gamma$ be a defect pointed group of $N_\beta$ such that $Q_\delta\leq P_\gamma\,;$ by
\cite[Proposition 5.3]{KP}, we have $Q=P\cap H$ and,  since we have
\cite[1.7]{KP}
$$\O G\,{\rm Tr}_N^G( b)\cong {\rm Ind}_N^G (\O N b)
\quad ,$$
it is easily checked that $P_\gamma$ is also a {\it defect pointed group\/} of $G_\alpha$ \cite[1.12]{P3}. Assume that the block $b$
is {\it nilpotent}; it follows from \cite[Proposition~6.5]{KP}  that
$b$ remains a {nilpotent block\/} of $H\.     R$ for any subgroup~$R$
of~$P\,,$ and from \cite[Theorem~6.6]{KP} that there is a unique
{local point\/} $\varepsilon$ of $R$ on $\O H$ such that~$R_\varepsilon\i P_\gamma\,.$

\bigskip \noindent {\bf £3.2.}\quad
Set ${\cal A} = \O Nb$ and ${\cal B} = \O Hb\,$. Choosing  $j\in \delta$
and $i\in \gamma$ such that $ij=ji=j\,,$ we set
\begin{center} ${\cal A}_\gamma=(\O G)_\gamma=i{\cal A} i$, ${\cal B}_\gamma=(\O H)_\gamma=i{\cal B}i$ and
${\cal B}_\delta=(\O H)_\delta=j{\cal B}j$. \end{center} Then ${\cal A}_\gamma$ is
a $P$-interior algebra with the group homomorphism $P\rightarrow
{\cal A}_\gamma^*$ mapping $u$ onto $ui$ for any $u\in P$, ${\cal
B}_\gamma$ is a $P$-stable subalgebra of ${\cal A}_{\gamma}$ and
${\cal B}_\delta$ is a $Q$-interior algebra with the group
homomorphism $Q\rightarrow {\cal B}_\delta^*$ mapping $v\in Q$ onto
$vj$ for any $v\in Q$. Clearly $\cal A$ is an $N/H$-graded algebra
with the $\bar x$-component $\O H x b$, where $\bar x\in N/H$ and
$x$ is a representative of $\bar x$ in $N$. Since $i$ belongs to the
$1$-component $\cal B$, ${\cal A}_\gamma$ is an $N/H$-graded algebra
with the $\bar x$-component $i(\O H x)i$.

\bigskip \noindent {\bf £3.3.}\quad  In \cite{KP} K\"ulshammer and Puig
describe the structure of any block of $G$ lying over $b$ in
terms of a new finite group $L$ which need not be {involved\/}
in $G$ \cite[Theorem~1.8]{KP}. More explicitly, $L$ is a {group
extension\/} of $\bar N$ by~$Q$ holding {strong uniqueness\/}
properties. In order to prove these properties, in \cite{KP} the
group $L$ is exhibited inside a suitable $\O\-$algebra
\cite[Theorem~8.13]{KP}, demanding a huge effort. But, as a matter
of fact, these properties can be obtained {directly\/} from the
so-called {local structure\/} of $G$ over $\O H b\,,$ a fact
that we only have understood recently. Then, with these  uniqueness
properties in hand, the second main result \cite[Theorem~1.12]{KP}
follows quite easily. With the notation and framework of \cite{KP},
we completely develop both new proofs.

\bigskip
\noindent {\bf £3.4.}\quad Denote by $\E_{(b,H,G)}$ the category ---
called the {\it extension category\/} associated to $G\,,$ $H$ and
$b$ --- where the objects are all the subgroups of $P$ and, for  any
pair of subgroups $R$ and $T$ of $P\,,$ the morphisms from $T$ to
$R$ are the pairs $(\psi_x,\bar x)$  formed by an injective group
homomorphism $\psi_x\,\colon T\to R$ and an element $\bar x$ of
$\bar N$ both determined by an element $x\in N$ fulfilling $T_\nu\i
(R_\varepsilon)^x$ where $\varepsilon$ and $\nu$ are the respective
local points of $R$ and~$T$ on $\O H$ determined by~$P_\gamma$ ---
in general, we should consider the {\it $(b,N)\-$Brauer pairs\/}
over the {$p\-$permutation $N\-$algebra $\O H b$\/}
[5,~Definition~1.6 and Theorem~1.8] but, in our situation, they
coincide with the  {local pointed groups\/} over this
$N\-$algebra. The composition in $\E_{(b,H,G)}$ is determined by the
composition of group homomorphisms and by the product in $\bar N \,.$

\bigskip
\noindent {\bf Theorem  £3.5.}\quad {\it There is a triple formed by
a finite group $L$ and by two  group homomorphisms
$$\tau : P\too L\quad and\quad \bar\pi : L\too \bar N
\leqno  £3.5.1\phantom{.}$$ such that $\tau$ is injective, that
$\bar\pi$ is surjective, that we have  ${\rm Ker}(\bar\pi) = \tau (Q)$ and
$\bar\pi\big(\tau (u)\big) =\bar u$ for any $u\in P\,,$ and
that these homomorphisms induce an equivalence of categories
$$\E_{(b,H,G)}     \cong \E_{(1,\tau (Q),L)}.
\leqno  £3.5.2$$

\smallskip
\noindent Moreover, for another such a triple $L'\,,$ $\tau'$ and
$\bar\pi'\,,$ there is a group isomorphism $\lambda\,\colon L\cong
L'\,,$ unique up to conjugation, fulfilling
$$\lambda\circ\tau = \tau'\quad and \quad \bar\pi'\circ\lambda = \bar\pi\,.$$\/}

\par
\noindent {\bf Proof:} Set $Z = Z(Q)\,,$  $M = N_G (Q_\delta)$ and
$\E = \E_{(b,H,G)}\,,$ denote by $\E (R,T)$ the set of
$\E\-$morphisms from $T$ to $R\,,$ and write $\E (R)$ instead of $\E
(R,R)\,;$ by the very definition of the category $\E\,,$ we have the
exact sequence
$$1\too C_H (Q)\too M\too \E (Q) \too 1 ;
$$ it is clear that $M$ contains $P$ and that we have
$C_H (Q)\cap P = Z\,;$ moreover, denoting by $\E_P (Q)$ the image of
$P$ in $\E (Q)\,,$ it is easily checked from
\cite[~Proposition~5.3]{KP} that $\E_P (Q)$ is a Sylow $p\-$subgroup
of~$\E (Q)\,.$

\smallskip
 We claim that the element $\bar h$ induced by $P$ in the second cohomology group ${\Bbb H}^2 \big(\E_P (Q),Z\big)$ belongs to the  image of~${\Bbb H}^2 (\E (Q),Z)\,.$ Indeed, according to in~
 \cite[~Ch.~XII,~Theorem~10.1]{CE},
 it suffices to prove that, for any subgroup $R$ of $P$ containing~$Z$ and any $(\varphi_x,\bar x)\in \E (Q)$ such that
$$(\varphi_x,\bar x)\circ\E_{\! R}(Q)\circ (\varphi_x,\bar x)^{-1}\i \E_{\!P} (Q),
\leqno  £3.5.3$$ the restriction ${\rm res}_{(\varphi_x,\bar x)}
(\bar h)$ of $\bar h$ {via\/} the conjugation by
$(\varphi_x,\bar x)$ and the element of ${\Bbb
H}^2\big(\E_R(Q), Z \big)$ determined by $R$ coincide; actually, we may
assume that $R$ contains $Q\,.$ Thus, $x$ normalizes $Q_\delta$ and
inclusion~£3.5.3 forces
$$C_H(Q)\.    R \i \big(C_H (Q)\.    P\big)^x;
$$ in particular, respectively denoting by $\lambda$
and $\mu$ the points of  $C_H (Q)\.    P$ and $C_H (Q)\.    R$ on $\O H$
such that  $\big(C_H (Q)\.    P\big)_\lambda$ and $\big(C_H
(Q)\.    R\big)_\mu$ contain $Q_\delta$ \cite[Lemma~3.9]{P4}, by
uniqueness we have
$$\big(C_H (Q)\.    R\big)_\mu\i \big(C_H (Q)\.    P\big)_\lambda
$$
 and, with the notation above, it follows from \cite[~Proposition~3.5]{KP} that
$P_\gamma$ and $R_\varepsilon$ are {defect pointed groups\/} of
the respective pointed groups $\big(C_H (Q)\.    P\big)_\lambda$ and
$\big(C_H (Q)\.    R\big)_\mu\,;$ consequently, there is $z\in C_H (Q)$
fulfilling $R_\varepsilon \i (P_\gamma)^{zx}$
\cite[~Theorem~1.2]{P5}. That is to say, the conjugation by $zx$
induces a group homomorphism $R\to P$ mapping $Z$ onto $Z$ and
inducing the element $(\psi_{zx},\overline{zx})$ of $\E (P,R)$ which
extends $(\varphi_x,\bar x)\,,$ so that the map
$${\rm res}_{(\varphi_x,\bar x)} : \Bbb H^2\big(\E_P (Q),Z\big)\too \Bbb H^2\big(\E_R (Q),Z\big)
$$ sends $\bar h$ to the element of~${\Bbb
H}^2\big(\E_R(Q), Z \big)$ determined by~$R$ \cite[~Chap.~XIV,
Theorem~4.2]{CE}.

\smallskip
 In particular, the corresponding element of ${\Bbb H}^2 \big(\E (Q),Z\big)$ determines a group extension
$$1\too Z\buildrel \tau \over\too L \buildrel \pi \over \longrightarrow \E (Q)\too 1
$$ and, since $\bar h\in \Bbb H^2\big(\E_P
(Q),Z\big)$ is the image of this element, there is  a {group
extension\/} homomorphism $\tau\,\colon P\to L$ \cite[~Chap.~XIV,
Theorem~4.2]{CE}; it is clear that $\tau$ is injective and, since
$\E_P (Q)$ is a Sylow $p$-subgroup of $\E (Q)\,,$ ${\rm Im}(\tau)$
is a Sylow $p\-$subgroup of $L\,;$ moreover, since $N = H\.     M$
\cite[Theorem~1.2]{P5}, we have
$$\bar N \cong M/C_H (Q)\.    Q\cong \E (Q)/\E_Q (Q)
;$$ in particular, $\pi$ determines a group homomorphism
$\bar\pi\,\colon L\to \bar N$ and, since $\tau$ is a {group
extension\/} homomorphism, we get $\bar\pi\big(\tau (u)\big) = \bar
u$ for any $u\in P$ and may choose $\pi$ in such a way that  we have
$$y\tau (v)y^{-1} = \tau\big(\varphi_x (v)\big)
\leqno  £3.5.4\phantom{.}$$ for any $y\in L$ and any $v\in Q$ where
$\pi (y) = (\varphi_x,\bar x)$ for some $x\in N\,.$ Then, we claim
that, up to a suitable modification of our choice of $\tau\,,$ the
group $L$ endowed  with $\tau$ and $\bar\pi$ fulfills the conditions
above; set $\hat\E = \E_{(1,\tau (Q),L)}$ for short.

\smallskip
For any pair of subgroups $R$ and $T$ of $P$ containing~$Q\,,$ since
we have $H\cap R = Q = H\cap T\,,$ denoting by $\varepsilon$ and $\nu$
the respective local points of~$R$ and $T$ such that $P_\gamma$ contains
$R_\varepsilon$ and $T_\nu\,,$ these local pointed groups contain~$Q_\delta$
and, in particular, any $\E\-$morphism
$$(\psi_x,\bar x) : T\too R
$$ determines an element $(\varphi_x,\bar
x)$ of $\E (Q)$ fulfilling
$$(\varphi_x,\bar x)\circ \E_T (Q)\circ (\varphi_x,\bar x)^{-1} \subset \E_R (Q)
\quad .$$
Thus, for any $y\in L$ such that $\pi (y) =
(\varphi_x,\bar x)\,,$ we have
$$y\,\tau (T)\,y^{-1}\i \tau (R)
\quad ;$$
more precisely, for any $w\in T$ and any $v\in Q\,,$ from equality~£3.5.4 we get
$$y\,\tau (v^w)\,y^{-1} = \tau \big(\varphi_x  (v^w)\big)
 = \tau \big(\varphi_x(v)\big)^{\tau (\psi_x (w))}
\quad ;$$
 moreover, since $x T x^{-1}\i R\,,$ we have
 $$\bar\pi \big(y\,\tau (w)\,y^{-1}\big) = \bar x\,\bar w\, \bar x^{-1}
 = \bar\pi \Big(\tau (\psi_x (w))\Big)
 \quad .$$
 Hence, for any $w\in T$ and a suitable $\theta_x (w)\in Z\,,$ we get
 $$y\,\tau \big(w\,\theta_x (w)\big)\,y^{-1} = \tau (\psi_x (w))
\quad .$$

\smallskip
 Conversely, since $R$ and $T$ have a unique (local) point on $\O Q\,,$ any
 $\hat\E\-$morphism from $T$ to $R$ induced by an element $y$ of $L$
 determines an element $\pi (y) = (\varphi_x,\bar x)$ of $\E (Q)\,,$ for some
 $x\in N\,,$ which still fulfills
 $$(\varphi_x,\bar x)\circ \E_T (Q)\circ (\varphi_x,\bar x)^{-1} \subset \E_R (Q)
\quad ;$$
thus, as above,  $x$ normalizes $Q_\delta$ and this inclusion
forces
$$C_H(Q)\.    T \i \big(C_H (Q)\.    R\big)^x
\quad .$$
Once again, respectively denoting by $\lambda$ and $\mu$ the
points of  $C_H (Q)\.    R$ and $C_H (Q)\.    T$ on $\O H$ such that
 $\big(C_H (Q)\.    R\big)_\lambda$ and $\big(C_H (Q)\.    T\big)_\mu$ contain
 $Q_\delta$  \cite[ Lemma~3.9]{P4}, and by $\varepsilon$ and $\nu$ the local points of $R$ and $T$
 on $\O H$ such that $P_\gamma$ contains $R_\varepsilon$ and~$T_\nu\,,$  it follows from
 \cite[~Proposition~3.5]{KP} that
 $R_\varepsilon$ and $T_\nu$ are defect pointed groups
of the respective pointed groups $\big(C_H (Q)\.    R\big)_\lambda$ and
$\big(C_H (Q)\.    T\big)_\mu\,;$ since by uniqueness we have
$$\big(C_H (Q)\.    T\big)_\mu\i \big(C_H (Q)\.    R\big)_\lambda
,$$
 there is $z\in C_H (Q)$ fulfilling  $T_\nu \i (R_\varepsilon)^{zx}$  \cite[Theorem~1.2]{P5}. That is to say, the conjugation by $zx$ induces a group homomorphism $\psi_{zx}\,\colon T\to R$ mapping $Z$ onto $Z$ and inducing the element $(\psi_{zx},\overline{zx})$ of $\E (R,T)$ which extends $(\varphi_x,\bar x)\,;$ hence, as above, for any $w\in T$ and a suitable $\theta_y (w)\in Z$ we get
 $$y\,\tau \big(w\,\theta_y (w)\big)\,y^{-1} = \tau (\psi_{zx} (w)).
 \leqno   £3.5.5$$

\smallskip
We claim that, for a suitable choice of $\tau\,,$ the elements
$\theta_x (w)$ and $\theta_y (w)$ are always trivial; then, the
equivalence of categories~£3.5.2 will be an easy consequence of the
above correspondences. Above, for any~$y\in L$ such that $\tau
(T)\subset \tau (R)^y$ we have found an element $\big(\psi_y,\bar\pi
(y)\big)\in \E (R,T)$ lifting~$\pi (y)$ in such a way that, for
any~$w\in T\,,$ we have
$$\tau \big(w\,\theta_y (w)\big) = \tau\big(\psi_y (w)\big)^y
\leqno   £3.5.6\phantom{.}$$ for a suitable $\theta_y (w)\in Z\,;$
note that, according to equality~£3.5.4, for any $v\in Q$ we have
$\theta_y (v) = 1\,,$ and whenever $y$ belongs to $\tau (R)$ we may
choose $\psi_y$ in such a way that  $\theta_y (w) = 1\,.$

\smallskip
In this situation, for any $w,w'\in T\,,$ we get
\begin{eqnarray*}
  \tau \big(ww'\theta_y (ww')\big) &=& \tau \big(\psi_y (ww')\big)^y \\
   &=& \tau\big(\psi_y (w) \big)^y\tau\big(\psi_y (w') \big)^y \\
   &=& \tau \big(w\,\theta_y (w)\big)\tau \big(w'\theta_y (w')\big) \\
   &=& \tau\big(w\,\theta_y (w)\,w'\theta_y (w')\big) \\
   &=& \tau \big(ww'\theta_y (w)^{w'}\theta_y (w')\big)
\end{eqnarray*}
and therefore, since $\tau$ is injective, we still get
$$\theta_y (ww') = \theta_y (w)^{w'}\theta_y (w')
\quad ;$$
in particular, for any $z\in Z$ we have
$$\theta_y (wz) = \theta_y (w)^z\,\theta_y (z) = \theta_y (w)
\quad .$$
In other words, the map $\theta_y$ determines a $Z\-$valued
$1\-${cocycle\/} from the image $\tilde T$ of~$T$ in
$\widetilde{\rm Aut}(Q) = {\rm Out} (Q)\,.$

\smallskip
 Actually, the {cohomology class\/}~$\bar\theta_y$ of this $1\-$cocycle does not depend on
 the choice of~$\psi_y\,;$ indeed, if another choice $\psi'_y$ determines~$\theta'_y\,\colon T\to Z$
 then we clearly have  $\psi'_y (T) = \psi_y (T)$ and, according to our argument above, there is $z\in C_H (Q)$ such that
$$(T_\nu)^z = T_\nu\qq \psi'_y = \psi_y\circ \kappa_z
\quad ,$$
where $\kappa_z\,\colon T\to T$ denotes the conjugation by
$z\,;$ actually, we still have
$$[z,T]\i H\cap T = Q
\quad .$$
But, since $T_\nu$ is a defect pointed group of $\big(C_H
(Q)\.    T\big)_\mu$ and, according to [4,~Theorem~1.2] and
\cite[Proposition~6.5]{KP}, $\mu$ determines a {nilpotent
block\/} of the group $C_H (Q)\.    T\,,$ we have $N_{C_H (Q)\.    T}
(T_\nu) = C_H (T)\.     T\,.$ Thus, $z$ belongs to~$Z\.    C_H (T)$ and we
actually may assume that $z$ belongs to $Z\,.$

\smallskip
 In this case, it follows from equality~£3.5.6 applied twice that
\begin{eqnarray*}
  \tau \big(w\,\theta'_y (w)\big) &=& \tau\big(\psi'_y (w)\big)^y \\
   &=& \tau\big(\psi_y (z w z^{-1})\big)^y \\
   &=& \tau \big((z w z^{-1})\, \theta_y (z w z^{-1})\big)
\end{eqnarray*}
for any $w\in T$ and, since $\theta_y (z w z^{-1}) =\theta_y (w)$ and $\tau$ is injective,  we get
$$\theta'_y (w)\theta_y (w)^{-1} = w^{-1}z w z^{-1} = (z^{-1})^w z
\quad .$$
 Consequently, denoting by $\T_{\!L}$ the category where the objects are the subgroups of $\tau (P)$
 and the set of morphisms $\T_{\! L}  \big(\tau (R),\tau (T)\big)$ from $\tau (T)$ to $\tau (R)$ is
 just the corresponding {\it transporter\/} in $L\,,$ the correspondence sending an element $y\in \T_{\! L}  \big(\tau (R),\tau (T)\big)$ to the cohomology class
 $\bar\theta_y$ of $\theta_y$  determines a map
$$ \bar\theta_{_{R,T}} : \T_{\! L} \big(\tau (R),\tau (T)\big)\too {\Bbb H}^1
(\tilde T,Z)
\quad .$$

\smallskip
Moreover, if $U$ is a subgroup of $P$ containing $Q$ and $t$ an
element of~$L$ fulfilling $\tau (U)\subset \tau (T)^t\,,$ as above
we can choose $\big(\psi_t,\bar\pi (t)\big)\in \E (T,U)$
lifting~$\pi (t)$ in such a way that, for any~$u\in U\,,$ we have
$$\tau \big(u\,\theta_t (u)\big) = \tau\big(\psi_t (u)\big)^t
$$ for a suitable $\theta_t (u)\in Z\,;$ then, the composition
$\big(\psi_y,\bar\pi (y)\big)\circ\big(\psi_t,\bar\pi (t)\big)$
lifts $\pi (yt)$ and, for any $u\in U\,,$ we may assume that
(cf.~£3.5.4)

\begin{eqnarray*}
  \tau \big(u\,\theta_{yt} (u)\big) &=& \tau\big((\psi_y\circ\psi_t) (u)\big)^{yt} \\
   &=& \tau \Big(\psi_t (u)\,\theta_y\big(\psi_t (u)\big)\Big)^t \\
   &=& \tau \big(u\,\theta_t (u)\big) \tau \Big(\theta_y
\big(\psi_t (u)\big)\Big)^t \\
   &=& \tau \bigg(u\,\theta_t (u)\,\pi (t)^{-1} \Big(\theta_y
\big(\psi_t (u)\big)\Big)\bigg)\quad ;
\end{eqnarray*}
finally, since $\tau$ is injective, using {additive notation\/}
in $Z$ we get
$$\theta_{yt} (u) = \theta_t (u) + \pi (t)^{-1}\Big(\theta_y
\big(\psi_t (u)\big)\Big)
\quad .$$

\smallskip
Hence, denoting by $\tilde t$ the image of $t$ in $\widetilde{\rm
Aut}(Q)$ and by $\psi_{\tilde t}\,\colon \tilde U\to \tilde T$ and
$\Z(\tilde t)\,\colon Z\cong Z$ the corresponding group
homomorphisms, we get the~{$1\-$cocycle condition\/}
$$\bar\theta_{yt}  = \bar\theta_t  + {\Bbb H}^1 \big(\psi_{\tilde t},
\Z (\tilde t)\big) (\bar\theta_y)
\quad ; \leqno  £3.5.7$$
in particular, since $\theta_y (w) = 0$ whenever $y\in\tau (R)\,,$ it is easily
checked from this condition that $\bar\theta_y$ only depends on the
class of $y$ in the {\it exterior quotient\/}
$$\tilde\T_{\! L} \big(\tau (R),\tau (T)\big) = \tau (R)\backslash \T_{\! L}
\big(\tau (R),\tau (T)\big) .$$ Thus, respectively denoting by
$\tilde L\,,$ $\tilde R\,,$ $\tilde T$ and $\tilde P$ the images of
$L\,,$ $\tau(R)\,,$ $\tau(T)$ and~$\tau (P)$ in $\widetilde{\rm
Aut}(Q)\,,$ the map $\bar\theta_{_{R,T}}$ above admits a
factorization
$$\skew4\tilde{\bar\theta}_{_{\tilde R,\tilde T}} : \tilde \T_{\!\tilde L}
(\tilde R,\tilde T)\too {\Bbb H}^1 \big(\tilde T,Z\big) .$$

\smallskip
That is to say, let us consider the {exterior quotient\/}
$\tilde\T_{\!\tilde L}$
 of the category $\T_{\!\tilde L}$ and the {contravariant\/} functor
$${\frak h^1_Z }: \tilde\T_{\!\tilde L}\too \Ab
$$
to the category of Abelian groups $\Ab$ mapping~$\tilde T$  on~${\Bbb H}^1
\big(\tilde T,Z\big)\,;$
then, identifying the $\tilde\T_{\!\tilde L}\-$morphism $\tilde y\in
\tilde \T_{\!\tilde L} (\tilde R,\tilde T)$ with the obvious {$\tilde\T_{\!\tilde L}\-$chain\/} $\Delta_1\too \tilde\T_{\tilde L}$
--- the {functor\/} from the category $\Delta_1\,,$ formed by
the objects $0$ and $1$ and a non-identity morphism $0\bullet 1$
from $0$ to $1\,,$ mapping $0$ on $T\,,$ $1$ on $R$ and $0\bullet 1$
on~$\tilde y$ --- the family $\bar \theta = \{\bar\theta_{\tilde
y}\}_{\tilde y}\,,$ where $\tilde y$ runs over the set of all the
$\tilde\T_{\!\tilde L}\-$morphisms, defines a {$1\-$cocycle\/}
from $\tilde \T_{\!\tilde L}$ to ${\frak h^1_Z}$ since
equalities~£3.5.7 guarantee that the {differential map\/} sends
$\bar\theta$ to zero.

\smallskip
We claim that this {$1\-$cocycle\/} is a {$1\-$coboundary\/}; indeed, for any subgroup~$\tilde R$ of $\tilde
P\,,$ choose a set of representatives $\tilde X_{\tilde R}\subset
\tilde L$ for the set of {double classes\/} $\tilde P\backslash
\tilde L/\tilde R$ and,
 for any $\tilde n\in  \tilde X_{\tilde R}\,,$ set $\tilde R_{\tilde n}
= \tilde R\cap P^{\tilde n}\,,$ consider the $\tilde\T_{\tilde
L}\-$morphisms $\tilde n\,\colon \tilde R_{\tilde n}\to \tilde P$
and $\tilde\imath_{\tilde R_{\tilde n}}^{\tilde R} \,\colon \tilde
R_{\tilde n} \to \tilde R$ respectively determined by $\tilde n$ and
by the trivial element of~$\tilde L\,,$ and denote by
$$({\frak h}^1_Z)^{^{\!\circ}}  (\tilde\imath_{\tilde R_{\tilde n}}^{\tilde R})
: {\Bbb H}^1 \big(\tilde R_{\tilde n},Z\big)\too {\Bbb H}^1
\big(\tilde R,Z\big) $$
 the corresponding {transfer homomorphism\/}  \cite[~Ch.~XII,~\S8]{CE}; then, setting
$$\bar\sigma_{\tilde R} = {\vert P\vert \over \vert L\vert}\.    \sum_{\tilde n\in
\tilde X_{\tilde R}} \big(({\frak h}^1_Z)^{^{\!\circ}}
(\tilde\imath_{\tilde R_{\tilde n}}^{\tilde R})\big)
(\bar\theta_{\tilde n})
\quad ,$$
 we claim that, for any ${\tilde y}\in\bar\T_{\!\tilde L}(\tilde R,\tilde T)\,,$ we
have
$$\bar\theta_{\tilde y} =
\bar\sigma_{\tilde T} - \big({\frak h^1_Z} (\tilde
y)\big)(\bar\sigma_{\tilde R})
\quad  . \leqno  £3.5.8$$

\smallskip
Indeed, note that ${\frak h^1_Z} (\tilde y)$ is the composition of the
restriction {via\/} the $\tilde\T_{\tilde L}\-$morphism
$$\tilde\imath_{\tilde y\tilde T \tilde y^{-1}}^{\tilde R} :
\tilde y\,\tilde T \,\tilde y^{-1}\too \tilde R $$ determined by the
trivial element of $L\,,$ with the conjugation determined by~$\tilde
y\,,$ which we denote by~${\frak h^1_Z} (\tilde y_*)\,;$
thus, by the corresponding {Mackey equalities\/}
 \cite[Ch.~XII,~Proposition~9.1]{CE}, we get
\begin{eqnarray*}
  &{\frak h^1_Z} (\tilde y)\Big(\sum_{\tilde n\in \tilde X_R}
\big(({\frak h}^1_Z)^{^{\!\circ}} (\tilde\imath_{\tilde R_{\tilde
n}}^{\tilde R})\big) (\bar\theta_{\tilde n})\Big) \\ =& {\frak
h^1_Z} (\tilde y_*)\Big(\sum_{\tilde n\in \tilde X_{\tilde R}}\,
\sum_{\tilde r\in \tilde Y_{\tilde n}}  \big(({\frak
h}^1_Z)^{^{\!\circ}} (\tilde\imath_{\tilde P^{\tilde n \tilde r}
\,\cap\, \tilde y\,\tilde T\, \tilde y^{-1}}^{\tilde y\, \tilde
T\,\tilde y^{-1}}) \circ {\frak h^1_Z} (\tilde r)\big)
(\bar\theta_{\tilde n})\Big) \\
   =& \sum_{\tilde n\in \tilde X_{\tilde R}}\, \sum_{\tilde r\in
\tilde Y_{\tilde n}} \big(({\frak h}^1_Z)^{^{\!\circ}}
(\tilde\imath_{\tilde P^{\tilde n\tilde r \tilde y} \,\cap\,\tilde
T}^{\tilde T}) \circ {\frak h^1_Z} (\tilde
r\tilde y)\big) (\bar\theta_{\tilde n})\quad ,
\end{eqnarray*}
where, for any $\tilde n\in X_{\tilde R}\,,$ the subset  $\tilde Y_{\tilde
n}\subset  \tilde R$ is a set of representatives for the set of {double classes\/} $\tilde R_{\tilde n} \backslash \tilde R/\,\tilde
y\,\tilde T\, \tilde y^{-1}$ and, for any $\tilde r\in \tilde
Y_{\tilde n}\,,$ we consider the $\tilde\T_{\tilde L}\-$morphisms
$$\tilde r: \tilde P^{\tilde n \tilde r} \cap \tilde y\,\tilde T\,\tilde y^{-1}
\too \tilde R_{\tilde n}\qq \tilde r\tilde y : \tilde P^{\tilde
n\tilde r \tilde y}\cap\tilde T \too \tilde R_{\tilde n}
\quad .$$

\smallskip
Moreover, setting $\tilde m = \tilde n\tilde r\tilde y$ for~$\tilde
n\in \tilde X_{\tilde R}$ and $\tilde r\in \tilde Y_{\tilde n}\,,$
since we assume that $\theta_{\tilde r} = 0\,,$ it follows from
equality~£3.5.7 that
$$ \big({\frak h^1_Z} (\tilde r\tilde y)\big) (\bar\theta_{\tilde n}) =
\bar\theta_{\tilde m} - \bar\theta_{\tilde r\tilde y}
=\bar\theta_{\tilde m} - \big({\frak h^1_Z}(\tilde\imath_{\tilde
T_{\tilde m}}^{\tilde T})\big)(\bar\theta_{\tilde y})
\quad ;$$
 thus, choosing $\tilde X_{\tilde T} = \bigsqcup_{\,\tilde n\in \tilde X_{\tilde R}} \tilde n\,\tilde Y_{\tilde n}\,\tilde y\,,$
 we get \cite[~Ch.~XII,~\S8.(6)]{CE}
\begin{eqnarray*}
  \bar\sigma_{\tilde T} - \big({\frak h^1_Z} (\tilde y)\big) (\bar\sigma_{\tilde R}) &=& {\vert P\vert \over  \vert L\vert}\.    \sum_{\tilde m\in \tilde
X_{\tilde T}} \big(({\frak h}^1_Z)^{^{\!\circ}} (\tilde\imath_{
\tilde T_{\tilde m}}^{\tilde T})\big)\Big(\bar\theta_{\tilde m} -
\big({\frak h^1_Z} (\tilde r\tilde y)\big) (\bar\theta_{\tilde
n})\Big) \\
   &=& {\vert P\vert \over \vert L\vert}\.     \sum_{\tilde m\in
\tilde X_{\tilde T}} \big(({\frak h}^1_Z)^{^{\!\circ}}
(\tilde\imath_{\tilde T_{\tilde m}}^{\tilde T})\big)\Big(\big({\frak
h^1_Z}(\tilde\imath_{\tilde T_{\tilde m}}^{\tilde T})\big)
(\bar\theta_{\tilde y})\Big) \\
   &=& \sum_{\tilde m\in \tilde X_{\tilde T}}{\vert\tilde
T/ \tilde T_{\tilde m}\vert \over\vert\tilde L/\tilde P\vert}\.
\bar\theta_{\tilde y}=\bar\theta_{\tilde y} \quad .
\end{eqnarray*}

\smallskip
In particular, for any subgroup $\tilde R$ of $\tilde P\,,$ we get
$$\bar \sigma_{\tilde R} = \big({\frak h^1_Z}(\tilde\imath_{\tilde R}^{\tilde P})\big)(\bar \sigma_{\tilde P})
$$
 and the element $\bar\sigma_{\tilde P}\in \Bbb H^1(\tilde P,Z)$
can be lifted to a $1\-$cocycle $\sigma_{\tilde P}\,\colon \tilde
P\to Z$ which determines a group automorphism $\sigma\,\colon P\cong
P$ mapping $u\in P$ on $u\,\sigma_{\tilde P}(\tilde u)$ where
$\tilde u$ denotes the image of $u$ in $\tilde P\,;$ moreover,
according to equality~£3.5.8, in \noindent £3.5.5 we may choose
$$\theta_y (w) = \sigma_{\tilde P}(\tilde w)\big(\pi (y)\big)^{-1}\Big(\sigma_{\tilde P}
\big(\widetilde{\psi_y (w)}\big)\Big)^{-1} .$$
 Hence, replacing $\tau$  by $\hat\tau = \tau\circ\sigma\,,$ the
maps $\pi$ and $\hat\tau$ still fulfill the conditons above and, for
any $w\in T\,,$ in~equality~£3.5.6 we get
\begin{eqnarray*}
  \tau\big(\psi_y (w)\big)^y &=& \tau \big(w\,\theta_y (w)\big) \\
  &=& \tau\bigg(w\big(w^{-1}\sigma (w)\big)\big(\pi (y)\big)^{-1}
\Big(\psi_y (w)^{-1} \sigma \big(\psi_y (w)\big)\Big)^{-1}\bigg) \\
   &=& \tau \bigg(\sigma (w)\big(\pi (y)\big)^{-1}\Big(\sigma
\big(\psi_y (w) \big)^{-1}\psi_y (w)\Big)\bigg) \\
   &=& \hat\tau (w)\tau \Big(\sigma\big(\psi_y
(w)\big)^{-1} \psi_y (w)\Big)^y \\
   &=& \hat\tau (w) \hat\tau\big(\psi_y
(w)^{-1}\big)^y \tau\big(\psi_y (w) \big)^y
\end{eqnarray*}
so that, as announced, we obtain
$$\hat\tau\big(\psi_y (w)\big)^y = \hat\tau (w)
\quad .$$

\smallskip
In conclusion, we get a functor from $\hat\E$ to $\E$  mapping any
$\hat \E\-$morphism
$$(\kappa_y,\bar y) : \hat\tau (T)\too \hat\tau (R)
$$
induced by an element $y$ of $L\,,$ where $\kappa_y$ denotes the corresponding conjugation by~$y$ which actually fulfills $\hat\tau (Q\.    T)\i
\big(\hat\tau (Q\.    R)\big)^y\,,$ on the $\E\-$morphism
$$\big(\psi_y,\bar\pi (y)\big) : T\too R
$$
 where $\psi_y\,\colon T\to R$ is the group homomorphism determined  by the equality
$$\hat\tau_R \circ\psi_y  = \kappa_y\circ \hat\tau_T
\quad ,$$
$\hat\tau_R$ and $\hat\tau_T$ denoting the respective
restrictions of $\hat\tau$ to $R$ and $T\,;$ indeed, it is clear
that this correspondence maps the composition of $\hat\E\-$morphisms
on the corresponding
 composition of $\E\-$morphisms. Moreover, it is clear that this functor is
 {faithful\/}, and it follows from our argument above that any $\E\-$morphism
 $$(\psi_x,\bar x) : T\too R
$$ comes from an $\hat\E\-$morphism from
$\hat\tau (T)$ to $\hat\tau (R)\,.$

\smallskip
Moreover, for another triple $L'\,,$ $\tau'$ and $\bar\pi'$
fulfilling the above conditions, the corresponding equivalences of
categories~£3.5.2 induce an equivalence of categories
$$\hat\E\cong \E_{(1,\tau' (Q),L')} = \E'
\quad ;\leqno  £3.5.9$$ in particular,  we have
a group homomorphism
$$\bar\sigma : L\too \hat\E \big(\hat\tau (Q)\big)\cong
\E' \big(\tau' (Q)\big)\cong L'/\tau' (Z) $$ and we claim that Lemma~£3.6
below applies to the finite groups $L$ and $L'\,,$ with the Sylow
$p\-$subgroup $\hat\tau (P)$ of $L\,,$ the Abelian normal $p\-$group
$\tau' (Z)$ of~$L'$ and the group homomorphism $\bar\sigma\,\colon L\to
L'/\tau' (Z)$ above; indeed, the group homomorphism $\hat\tau (P)\to L'$
mapping $\hat\tau (u)$ on $\tau' (u)\,,$ for any $u\in P\,,$ clearly
lifts the restriction of $\bar\sigma$ and it is easily checked from
the equi-valence~£3.5.9 that it fulfills condition~£3.6.1 below.
Consequently, the last statement immediately follows from this
lemma. We are done.

\bigskip
\noindent {\bf Lemma~£3.6.}\quad {\it Let $L$ be a finite group, $M$
a group, $Z$ a normal Abelian $p'\-$divisible subgroup of $M$~and
$\bar\sigma\,\colon L\to \bar M = M/Z$ a group homomorphism. Assume
that, for  a Sylow $p\-$subgroup   $P$ of $L\,,$ there exists a
group homomorphism $\tau\,\colon P\to M$ lifting the restriction of
$\bar\sigma$ to $P$ and fulfilling the following condition

\smallskip
\noindent {\bf £3.6.1}\quad For any subgroup $R$ of $P$ and any
$x\in L$ such that $R^x\subset P\,,$ there is $y\in M$ such that
$\bar\sigma (x) = \bar y$ and $\tau (u^x) = \tau (u)^y$~for any
$u\in R\,.$

\smallskip
\noindent Then, there is a group homomorphism $\sigma\,\colon L\to
M$ lifting $\bar\sigma$ and extending~$\tau\,.$ Moreover, if
$\sigma'\,\colon L\to M$ is a group homomorphism which lifts
$\bar\sigma$ and extends~$\tau\,,$ then there is $z\in Z$ such that
$\sigma' (x) = \sigma (x)^z$ for any $x\in L\,.$\/}

\bigskip
\noindent {\bf Proof:}  It is clear that $\bar\sigma$ determines an
action of $L$ on $Z$ and it makes sense to consider the {cohomology groups\/} $\Bbb H^n (L,Z)$ and $\Bbb H^n (P,Z)$ for any
$n$ in~$\Bbb N\,.$ But, $M$ determines an element $\bar\mu$
of~${\Bbb H}^2 (\bar M,Z)$  \cite[~Chap.~XIV, Theorem~4.2]{CE}
 and if  there is a group homomorphism $\tau\,\colon P\to M$ lifting the restriction of~$\bar\sigma$
 then the corresponding image of $\bar\mu$ in
${\Bbb H}^2 (P,Z)$ has to be zero  \cite[Chap.~XIV,
Theorem~4.2]{CE}; thus, since the restriction map
$${\Bbb H}^2 (L,Z)\too  {\Bbb H}^2 (P,Z)
$$
 is injective \cite[~Ch.~XII,~Theorem~10.1]{CE}, we also get
$$\big({\Bbb H}^2 (\bar\sigma,{\rm id}_Z)\big)(\bar\mu) = 0
$$
 and therefore there is a group homomorphism $\sigma\,\colon L\to M$ lifting $\bar\sigma\,.$

\smallskip
 At this point, the {difference\/} between $\tau$ and the restriction of $\sigma$
to $P$ defines a {$1\-$cocycle\/} $\theta\,\colon P\to Z$ and,
for any subgroup $R$ of $P$ and any $x\in L$ such that $R^x\subset
P\,,$ it follows from condition~£3.6.1 that, for a suitable $y\in M$
fulfilling $\bar y =  \bar\sigma (x)\,,$ for any $u\in R$ we have

\begin{eqnarray*}
  \theta (u^x) &=& \tau (u^x)^{-1}\sigma (u^x)\\
  & =& \tau (u^{-1})^y\sigma (u)^{\sigma (x)} \\
   &=& \tau (u^{-1})^y \tau (u)^{\sigma (x)}\theta (u)^{\sigma (x)} \\
   &=& \Big(\big(y\sigma
(x)^{-1}\big)^{-1}\big(y\sigma (x)^{-1}\big)^{\tau (u)}\theta
(u)\Big)^{\sigma (x)}\quad ;
  \end{eqnarray*}
consequently, since the map sending $u\in R$ to
$$\big(y\sigma (x)^{-1}\big)^{-1}\big(y\sigma (x)^{-1}\big)^{\tau (u)}\in Z
$$ is a {$1\-$coboundary\/},  the
cohomology class $\bar\theta$ of $\theta$ is $L\-${stable\/},
and it follows again from  \cite[~Ch.~XII,~Theorem~10.1]{CE} that it
is the restriction of a suitable element $\bar \eta\in {\Bbb H}^1
(L,Z)\,;$ then, it suffices to modify $\sigma$ by a representative
of~$\bar\eta$ to~get a new group homomorphism $\sigma'\,\colon L\to
M$ lifting $\bar\sigma$ and extending~$\tau\,.$

\smallskip
Now, if $\sigma'\,\colon L\to M$ is a group homomorphism which lifts
$\bar\sigma$ and extends~$\tau\,,$ the element $\sigma' (x)\sigma
(x)^{-1}$ belongs to $Z$ for any $x\in L$ and thus, we get a {$1\-$cocycle\/} $\lambda\, \,\colon L\to Z$ mapping $x\in L$ on
$\sigma' (x)\sigma (x)^{-1}\,,$ which vanish over~$P\,;$ hence, it
is a {$1\-$coboundary\/} \cite[~Ch.~XII,~Theorem~10.1]{CE} and
therefore there exists $z\in Z$ such that
$$\lambda (x) = z^{-1}\sigma (x) z\sigma (x)^{-1}
$$ so that we have $\sigma' (x) =\sigma (x)^z$ for any $x\in L\,.$ We are done.

\bigskip
\noindent {\bf £3.7.}\quad Since $Q$ normalizes a unitary {full
matrix\/} $\O\-$subalgebra $T$ of ${\cal B}_\delta$ such that
\cite[~Theorem~1.6]{P4}
$${\cal B}_\delta\cong T\,Q\qq {\rm rank}_\O (T)\equiv 1 \bmod p
\quad ,\leqno  £3.7.1$$
the action of $Q$ on $T$ admits a unique lifting to a group homomorphism
\cite[1.8]{P4}
$$Q\too {\rm Ker}({\rm det}_T)
\quad ;$$
hence, we have
$$\B_\delta\cong T\otimes_\O \O Q$$
and therefore ${\B}_\delta$ admits a unique two-sided ideal
${\frak n}_\delta$ such that, considering ${\cal B}_\delta/{\frak n}_\delta$ as a $Q$-interior $\O$-algebra, there is an isomorphism
$${\cal B}_\delta/{\frak n}_\delta\cong T$$ of $Q$-interior $\O$-algebras. Then, a canonical {embedding\/} $f_\delta\,\colon {\cal
B}_\delta\to {\rm Res}_Q^H ({\cal B})$ \cite[~2.8]{P4} and the ideal
${\frak n}_\delta$ determine a two-sided ideal ${\frak n}$ of $\cal B$
such that $S = {\cal B}/{\frak n}$ is also a {full matrix\/}
$\O\-$algebra.

\bigskip
\noindent {\bf Proposition~£3.8.} {\it With the notation above, the
action of $N$ on $\cal B$ stabilizes~${\frak n}\,.$\/}

\bigskip
\noindent {\bf Proof:} Since we have $N = H\.    N_G (Q_\delta)\,,$ for
the first statement we may consider $x\in N_G (Q_\delta)\,;$ then,
denoting by $\sigma_x$ the automorphism of $Q$ induced by the
conjugation by $x\,,$ it is clear that the isomorphism
$$f_x : {\rm Res}_{\sigma_x}\big({\rm Res}_Q^H ({\cal B})\big)\cong {\rm Res}_Q^H ({\cal B})
$$
of $Q$-interior algebras mapping $a\in \cal B$ on $a^x$ induces a commutative
diagram of {\it exterior\/} homomorphisms of  $Q$-interior
algebras \cite[2.8]{P4}
$$\matrix{{\rm Res}_{\sigma_x}\big({\rm Res}_Q^H ({\cal B})\big)&\buildrel \tilde f_x\over\cong &{\rm Res}_Q^H ({\cal B})\cr
\hskip-10pt{\scriptstyle \tilde
f_\delta}\hskip4pt\uparrow&\phantom{\Big\uparrow}&\uparrow\hskip4pt
{\scriptstyle \tilde f_\delta}\hskip-10pt\cr {\rm Res}_{\sigma_x}
({\cal B}_\delta)&\buildrel (\tilde f_x)_\delta\over\cong& {\cal
B}_\delta\cr}
\quad  ;$$ moreover, the uniqueness of ${\frak n}_\delta$
clearly implies that this ideal is stabilized by~$(\tilde f_x)_\delta\,;$
consequently, ${\frak n} $ is still stabilized by~$\tilde f_x\,.$

\bigskip
\noindent {\bf £3.9.} In particular, $N$ acts on the {full matrix\/}
$\O\-$algebra $S$ and therefore the action on $S$ of any element $x\in N$
can be lifted to a suitable $s_x\in S^*\,;$ thus, setting ${\rm r }=
{\rm rank}_\O(S)\,,$ denoting by $\bar H$ the image of $H$ in $S^*$ and considering a finite extension $\O'$ of $\O$ containing the group $U$ of
$\vert H\vert\-$th roots of unity and the ${\rm r}\-$th roots of ${\rm det}_S (s_x)$ for any $x\in N\,,$ since ${\rm r}$ divides $\vert H\vert\,,$ the {\it pull-back\/}
$$\matrix{N &\too & {\rm Aut}(\O'\otimes_\O S)\cr
\uparrow&\phantom{\big\uparrow}&\uparrow\cr
\hat N &\too &(U\otimes\bar H)\.{\rm Ker}({\rm det}_{\O'\otimes_\O S})\cr}$$
determines a central extension $\hat N$ of $N$ by $U\,,$ which clearly does not depend on the choice of $\O'\,;$ moreover, the inclusion $H\i N$ and the structural group homomorphism $H\to (\O'\otimes_\O S)^*$ induces
an injective group homomorphism $H\to \hat N$ with an image which
is a normal subgroup of~$\check N$ and has a {trivial\/}
intersection with the image of $U$
--- we identify this image with $H$ and set
$$\skew3\hat{\bar N} = \hat N/H\quad .$$
We will consider the $H\-$interior $N\-$algebras (see \cite[2.1]{P7})
$$\hat{\cal A} = S^\circ\otimes_\O {\cal A}\qq \hat{\cal B} =
S^\circ\otimes_\O {\cal B}$$
and note that $\O'\otimes \hat \A$ actually has an $\hat N\-$interior algebra structure.

\bigskip
\noindent {\bf £3.10.}\quad  On the other hand, since $b$ is also a
{nilpotent\/} block of the group $H\.    P\,,$ it is easily checked
that \cite[1.9]{P4}
$$\O(H\.    P)b\big/J\big(\O(H\.    P)b\big)\cong k\otimes_\O S
\quad ;$$
moreover, since the inclusion map $\O H\to \O (H\.    P)$ is a {semicovering of $P\-$algebras\/} \cite[Example~3.9, 3.10
and~Theorem~3.16]{KP}, we can identify $\gamma$ with a local point
of $P$ on $\O(H\.    P)b$. Set $\O(H\.    P)_\gamma=i(\O(H\.    P))i$ and
$S_\gamma=\bar\imath S\bar\imath$, where $\bar\imath$ is the image of $i$ in
$S\,;$ then, as in £3.7 above,  we have an isomorphism of $P$-interior algebras \cite[Theorem~1.6]{P4}
$$\O(H\.    P)_\gamma\cong S_\gamma\, P
\quad ,\leqno  £3.10.1$$
$S_\gamma$ is actually a {\it Dade $P\-$algebra\/}
--- namely, a {full matrix\/} $P\-$algebra over $\O$ where $P$
stabilizes an $\O\-$basis containing the unity element --- such that ${\rm rank}_\O(S_\gamma)\equiv 1\,\, {\rm mod}\,\, p$, and the action of $P$ on $S_\gamma$ can be uniquely lifted to a group
homomorphism $P\to {\rm
Ker}({\rm det}_{S_\gamma})$ \cite[1.8]{P4}, so that
isomorphism~£3.10.1 becomes
$$\O(H\.    P)_\gamma\cong S_\gamma\otimes_\O \O P¡£
\quad .\leqno  £3.10.2$$

\bigskip
\noindent {\bf Proposition~£3.11.}\quad {\it With the notation
above, the structural homomorphism ${\cal B}_\gamma\to S_\gamma$ of
$P\-$algebras is a strict semicovering.\/}

\bigskip
\noindent {\bf Proof:} It follows from isomorphism~£3.10.2 that the
canonical homomorphism of $P\-$algebras
$$\O(H\.    P)_\gamma\too S_\gamma
\leqno  £3.11.1\phantom{.}$$ admits a $P\-$algebra section mapping
$s\in S_\gamma$ on the image of $s\otimes 1$ by the inverse of that
isomorphism, which proves that the $P$-interior algebra
homomorphism~£3.11.1 is a {covering\/} \cite[4.14 and
Example~4.25]{P4}; thus, since the inclusion map $\O H\to \O (H\.    P)$
is a semicovering of $P\-$algebras, the canonical homomorphism of
$P\-$algebras
$${\cal B}_\gamma = (\O H)_\gamma\too S_\gamma
$$ remains a {semicovering\/} \cite[~Proposition~3.13]{KP}; moreover,
since ${\frak n}\i J(\cal B)\,,$ it is a {strict semicovering\/}
\cite[~3.10]{KP}.

\bigskip
\noindent {\bf £3.12.}\quad  Consequently, it easily follows from
\cite[~Theorem~3.16]{KP} and \cite[~Proposition~5.6]{P4} that we
still have a {strict semicovering\/} homomorphism of
$P\-$algebras
$$(S_\gamma)^\circ\otimes_\O {\cal B}_\gamma\too (S_\gamma)^\circ\otimes_\O S_\gamma\cong {\rm End}_\O (S_\gamma)
\quad ;\leqno  £3.12.1$$
hence, denoting by $\hat \gamma$ the local point of $P$ over
$(S_\gamma)^\circ\otimes_\O {\cal B}_\gamma$
determined by~$\gamma\,,$ the image of $\hat\gamma$ in $(S_\gamma)^\circ\otimes_\O S_\gamma$ is contained in the corresponding local point of $P$
and therefore we get  a  {strict semicovering\/} homomorphism  \cite[~5.7]{P4}
$$\hat{\cal B}_{\hat\gamma}\too \O \cong ((S_\gamma)^\circ\otimes_\O S_\gamma)_{\hat\gamma}$$
 of $P\-$algebras; that is to say, any $\hat\imath\in \hat \gamma$ is
actually a primitive idempotent
 in~$\hat{\cal B}$ and therefore, for any local pointed group $R_{\hat\varepsilon}$ over $\hat{\cal B}$ contained in~$P_{\hat\gamma}\,,$
it also belongs to $\hat\varepsilon\,;$ in particular,  denoting by
$\hat \delta$ the local point of $Q$ over
 $(S_\gamma)^\circ\otimes_\O  {\cal B}_\gamma$ determined by $\delta\,,$ we clearly have
 $\hat{\cal B}_{\hat\delta} = \hat\imath\hat{\cal B}\hat\imath\cong \O Q$ (cf.~£3.7.1).

\bigskip
\noindent {\bf £3.13.}\quad
As~in~\cite[~2.11]{KP}, we consider the $P$-interior algebra
$\hat{\cal A}_{\hat\gamma} = \hat\imath\hat{\cal A}\hat\imath\,;$ since $\cal A$ is an $N/H$-graded algebra, $\hat{\cal A}_{\hat\gamma}$ is also an $N/H$-graded algebra. On the other hand, since $\O'/J(\O')\cong k\,,$ we get a group
homomorphism $\varpi\,\colon U\to k^*$ and, setting $\Delta_\varpi (U)
= \{(\varpi (\xi),\xi^{-1})\}_{\xi\in U}\,,$ we obtain the obvious $k^*\-$group
$$\skew3\hat{\bar N}^{^k} =( k^*\times \skew3\hat{\bar N})/\Delta_\varpi (U)
\quad ; $$
then, with the notation of Theorem~£3.5, we  set \cite[~5.7]{P6}
$$\hat L = {\rm Res}_{\bar\pi} (\skew3\hat{\bar N}^{^k})
\quad ;\leqno £3.13.1$$
thus, $\O_*\hat L^{^\circ}$ becomes a
$P$-interior algebra {via\/} the lifting $\hat\tau\,\colon P\to \hat L^{^\circ}$ of the group
homomorphism $\tau\,\colon P\to L\,,$ and it has an
obvious $L/\tau(Q)$-graded algebra structure. The group homomorphism $\bar
\pi$ induces a group isomorphism $L/\tau(Q)\cong N/H$, through which
we identify $L/\tau(Q)$ and $N/H\,,$ so that $\O_*\hat L^{^\circ}$ becomes an
$N/H$-graded algebra.

\bigskip
\noindent {\bf Theorem~£3.14.}\quad {\it With the notation above, we
have a $P$-interior  and $N/H$-graded algebra isomorphism $\hat{\cal A}_{\hat\gamma}\cong \O_*\hat L^{^\circ}\,.$ }

\bigskip
\noindent {\bf Proof:} Choosing $\hat\imath\in \hat\gamma\,,$ we
consider the groups
$$M = N_{(\hat\imath\hat{\cal A} \hat\imath)^*} (Q\.    \hat\imath)/k^*\.    \hat\imath\qq Z
=  \big((\hat\imath\hat{\cal B}\hat\imath)^Q\big)^*\!
\big/k^*\.    \hat\imath\cong 1 + J\big(Z (\O Q)\big)
\quad ;$$
it is clear
that $Z$ is a normal Abelian $p'\-$divisible subgroup of $M\,,$ and
we set~$\bar M = M/Z\,.$ In order to apply Lemma~£3.6, let $R$ be a
subgroup of $P$ and $y$ an element of $L$ such that $\tau (R)\i \tau
(P)^y\,;$ since $\tau (Q)$ is normal in~$L\,,$ we actually may
assume that $R$ contains $Q\,.$ According to the equivalence of
categories~£3.5.2, denoting by $\varepsilon$ the unique local point
of $R$ on $\cal B$ fulfilling $R_\varepsilon\i P_\gamma$
\cite[~Theorem~6.6]{KP}, there is $x_y\in N$ such that
$$\bar x_y = \bar\pi (y)\quad,\quad R_\varepsilon\i (P_\gamma)^{x_y} \qq \tau ({}^{x_y}v) ={}^y\tau (v) \hbox{\ \  for any $v\in R$}
\quad ;\leqno  £3.14.1$$
in particular, $x_y$ normalizes $Q_\delta\,.$

\smallskip
By Proposition 3.11, a local pointed group $R_\varepsilon$ on $\cal B$ such that
$$Q_\delta\leq R_\varepsilon\leq P_\gamma$$
 determines a local pointed group $R_{\tilde \varepsilon}$ on $S$ through the composition
 $${\cal B}_\gamma\too S_\gamma\hookrightarrow S$$
  (see \cite[Proposition 3.15]{KP}). Since $S_\gamma$ has a $P$-stable $\O$-basis, $S_\varepsilon$ still has a $R$-stable $\O$-basis and, by \cite[Theorem 5.3]{P4}, there are  unique local pointed groups $R_{\tilde\varepsilon}$ on $S_\varepsilon$
  and $R_{\hat\varepsilon}$ on $\hat{\cal B}$ such that $\hat l(\tilde l\otimes l)=\hat l=(\tilde l\otimes l)\hat l$ for suitable $l\in \varepsilon$,
  $\tilde l\in \tilde\varepsilon$ and $\hat l\in \hat\varepsilon\,;$.
then, we claim that $R_{\hat\varepsilon}\i (P_{\hat\gamma})^{x_y}$ and
that $x_y$ stabilizes~$Q_{\hat\delta}\,.$ Indeed, since $(R_\varepsilon)^{x_y^{-1}}\i P_\gamma$, we have $(R_{\tilde\varepsilon})^{x_y^{-1}}\i P_{\tilde\gamma}$ and then it follows from \cite[Proposition 5.6]{P4} that we have
$(R_{\hat\varepsilon})^{x_y^{-1}}\i  P_{\hat\gamma}$ or, equivalently,
$R_{\hat\varepsilon}\i (P_{\hat\gamma})^{x_y}\,;$ moreover, since $\delta$ is the unique local point of $Q$ such that $Q_\delta$ is contained in~$P_\gamma\,,$ again by \cite[Proposition 5.6]{P4} we can easily conclude that $x_y$ stabilizes
$Q_{\hat\delta}\,.$

\smallskip In particular, since the image of $\hat\imath^{\,x_y}$ in $\hat{\cal B} (R_{\hat\varepsilon})$ is not zero [14,~2.7] and since $\hat\imath$ is primitive in
$\hat{\cal B}\,,$ $\hat\imath^{\,x_y}$ belongs to $\hat\varepsilon$ and therefore, since $\hat\imath$ also belongs to $\hat\varepsilon\,,$ there is
$\hat a_y\in (\hat{\cal B}^R)^*$ such that
$\hat\imath^{\,x_y} = \hat\imath^{\,\hat a_y}\,;$ choose $s_y\in S^*$ lifting
the action of $x_y$ on $S$ and set $\hat x_y = s_y\otimes x_y\,,$ so that we have
$$\hat\imath^{\,x_y} = (\hat x_y)^{-1} \hat\imath\, \hat x_y
\quad ;$$
then, since $\hat x_y$ and $\hat a_y$ normalize~$Q\,,$
the element $\hat x_y    \hat a_y^{-1}$ of $\hat A$ normalizes $Q\. \hat\imath$ and therefore $\hat x_y    \hat a_y^{-1}\hat\imath$
determines an element $m_y$ of~$M\,.$ We claim that the image $\bar m_y$ of $m_y$ in $\bar M$ only depends
on~$y\in L$ and that, in the case where $R_\varepsilon =Q_\delta\,,$ this correspondence determines a group homomorphism
$$\bar\sigma : L\too \bar M
\quad .$$

\smallskip
Indeed, if $x'\in N$ still fulfills conditions~£3.14.1 then we necessarily have
$x' = x_y\,z$ for some $z\in C_H (R)$ and therefore it suffices to choose
the element $\hat a_y\.    z$ of $(\hat B^R)^*$ in the definition above.
On the other hand, if  $\hat a'\in (\hat B^R)^*$ still fulfills
$\hat\imath^{\,\hat x_y} = \hat\imath^{\,\hat a'}$ then we clearly have
$\hat a' = \hat c\,\hat a_y$ for some $\hat c\in (\hat B^R)^*$ centralizing
$\hat\imath\,,$ so that $\hat c\,\hat\imath$ belongs to
$(\hat\imath\hat B\hat\imath)^Q\,;$ hence, the image of $\hat x_y\hat a_y^{-1}\hat c^{-1}\hat\imath$ in $\bar M$
coincides with $\bar m_y\,.$ Moreover, in the case where $R_\varepsilon
=Q_\delta\,,$ for any element $y'$ in $L$
we clearly can choose $\hat x_{yy'} = \hat x_y\, \hat x_{y'}\,;$ then, we have
$$\hat\imath^{\,\hat x_{yy'}} = (\hat\imath^{\,\hat a_y})^{\hat x_{y'}} = \hat\imath^{\hat x_{y'}\.    (\hat a_y)^{\hat x_{y'}}}
= \hat\imath^{\hat a_{y'}(\hat a_y)^{\hat x_{y'}}}$$
and therefore, since $\hat a_{y'}(\hat a_y)^{\hat x_{y'}}$ still belongs to
$(\hat B^Q)^*\,,$ we clearly can choose
$\hat a_{yy'} = \hat a_{y'}(\hat a_y)^{\hat x_{y'}}\,,$ so that we get
$$\hat x_{yy'}\.     \hat a_{yy'}^{-1}\hat\imath = \hat x_y\,\hat x_{y'}\.    \big(\hat a_{y'}
(\hat a_y)^{\hat x_{y'}}\big)^{-1}\hat\imath = (\hat x_y\.    \hat a_y^{-1}\hat\imath)(\hat x_{y'}\.    \hat a_{y'}^{-1}\hat\imath)$$
which implies that $\bar m_{yy'} = \bar m_y\,\bar m_{y'}\,.$ This proves our claim.

\smallskip
In particular, for any $u\in P\,,$ we can choose $x_{\tau (u)} = u$
and $\hat a_{\tau (u)} = 1\,;$ moreover, since the action of $P$ on
$S_\gamma$ can be lifted to a unique group~homomorphism
$\varrho \,\colon P\to {\rm Ker}({\rm det}_{S_\gamma})$ \cite[~1.8]{P4},
we may choose $\hat x_{\tau (u)} = \varrho (u)\otimes u\,;$ then, it is clear that the correspondence $\tau^*$ mapping
$\tau (u)$ on the image of $(\varrho (u)\otimes u)    \hat\imath$ in $M$ defines a group
homomorphism from $\tau (P)\i L$ to~$M$ lifting the corresponding
restriction of $\bar\sigma\,.$

\smallskip
Finally, we claim that $\tau^*$ fulfills condition~£3.6.1; indeed,
coming back to the general inclusion $\tau (R)\i \tau (P)^y$ above,
we clearly have $\bar\sigma (y) = \bar m_y$ and, according to the
right-hand equalities
 in~£3.14.1, for any $v\in R$ we get
$$\tau^*\big(\tau (v)^y\big) = v^{x_y}\.    \hat\imath = (v\.    \hat\imath)^{m_y} =\tau^*\big(\tau (v)\big)^{m_y}
\quad .$$
Consequently, it follows from Lemma~£3.6 that $\bar\sigma$ can
be lifted to a group homomorphism $\sigma\,\colon L\to M$ extending
$\tau^*\,;$ moreover, the inverse image of~$\sigma (L)$ in
$N_{(\hat\imath\hat{\cal A} \hat\imath)^*} (Q\.\hat\imath)$ is a $k^*\-$group
which is clearly contained in
$$\hat N\.(\O'^*\otimes 1)\subset \O'\otimes_\O \hat\A
\quad ;$$
hence, according to definition~£3.13.1, $\sigma$ still can be lifted to a $k^*\-$group
homomorphism
$$\hat\sigma : \hat L^{^\circ}\too N_{(\hat\imath\hat{\cal A} \hat\imath)^*} (Q\.    \hat\imath)
$$ mapping $\tau (u)$ on $u\.    \hat\imath$
for any $u\in P\,;$ hence, we get a $P$-interior and $N/H$-graded algebra homomorphism
$$\O_*\hat L^{^\circ}\too \hat{\cal A}_{\hat\gamma}
\quad .\leqno  £3.14.2$$
We claim that homomorphism £3.14.2 is an
isomorphism.

\smallskip
Indeed, denoting by $X\i  N_G (Q_\delta)$ a set of representatives
for $\bar N = N/H\,,$ it is clear that we have
$${\cal A} = \bigoplus_{x\in X} x\.    {\cal B}
$$
and therefore we still have
$$\hat{\cal A} = S\otimes_\O {\cal A} = \bigoplus_{x\in X} (s_x\otimes x)     (S\otimes_\O {\cal B})  = \bigoplus_{x\in X} (s_x\otimes x)     \hat{\cal B}
\quad ;$$
moreover, choosing as above an
element $\hat a_x\in (\hat{\cal B}^Q)^*$ such that
$\hat\imath^{\, x} = \hat\imath^{\,\hat a_x}\,,$ it is clear
that $(s_x\otimes x)    \hat a_x^{-1}\hat{\cal B}  =(s_x\otimes x)   \hat{\cal B} $
for any $x\in X$ and therefore we get
$$\hat{\cal A}_{\hat\gamma} = \hat\imath \hat{\cal A}\hat\imath
= \bigoplus_{x\in X} ((s_x\otimes x)  \hat a_x^{-1}\hat\imath)
(\hat\imath\hat{\cal B}\hat\imath)
\quad ;$$ thus, since we know that $\hat\imath\hat{\cal B}\hat\imath\cong
\O Q$ and that $L/\tau (Q)\cong \bar N\,,$ denoting by $Y\i L$ a set
of representatives for $L/\tau (Q)$ and by $\hat y$ a lifting of
$y\in Y$ to $\hat L\,,$ we still get
$$\hat{\cal A}_{\hat\gamma} \cong \bigoplus_{y\in Y} \hat\sigma (\hat y)\, \O Q
$$ which proves that homomorphism~£3.14.2
is an isomorphism.

\bigskip
\noindent {\bf Corollary~£3.15.}\quad {\it With the notation above,
we have a  $P$-interior  and $N/H$-graded algebra isomorphism ${\cal A}_\gamma\cong S_\gamma\otimes_\O
\O_*\hat L^{^\circ}\,.$\/}

\bigskip
\noindent {\bf Proof:} Since $\hat{\cal A} = S^\circ\otimes_\O {\cal
A}$ and we have a  $P$-interior algebra embedding $\O\to S_\gamma\otimes_\O
S_\gamma^\circ$ \cite[~5.7]{P4}, we
still have the following commutative diagram of {\it exterior\/} $P$-interior algebra embeddings and homomorphisms
\cite[2.10]{KP}
$$\matrix{&{\cal A}_\gamma&\too&\hskip-20pt S_\gamma\otimes_\O S_\gamma^\circ \otimes_\O {\cal A}_\gamma&\cr
&\hskip-20pt\nearrow&\nearrow\hskip-60pt&\uparrow\cr {\cal
B}_\gamma&\too& S_\gamma\otimes_\O S_\gamma^\circ \otimes_\O {\cal
B}_\gamma& S_\gamma\otimes_\O \hat{\cal
A}_{\hat\gamma}&\hskip-20pt\cong&\hskip-20pt S_\gamma\otimes_\O
\O_*\hat L^{^\circ}\cr &&&\hskip-40pt\nearrow&\nearrow\hskip-20pt&\cr
&&S_\gamma\otimes_\O \hat{\cal
B}_{\hat\gamma}\hskip-40pt&\hskip-20pt\cong&\hskip-20pt
S_\gamma\otimes_\O \O Q\cr}
\quad ;\leqno  £3.15.1$$
moreover, since the unity element is primitive in $(S_\gamma)^P$ and the kernel of the
canonical homomorphism
$$(S_\gamma\otimes_\O \O Q)^P\too (S_\gamma)^P
$$
is contained in the radical,  the unity element is  primitive in $(S_\gamma\otimes_\O \O Q)^P$ too; since
$P$ has a unique local point over $ S_\gamma\otimes_\O
S_\gamma^\circ \otimes_\O {\cal A}_\gamma$
\cite[Proposition~5.6]{P4}, from diagram~£3.15.1 we get the
announced isomorphism.

\bigskip
\noindent {\bf £3.16.}\quad Let us take advantage of this revision
to correct the erroneous  proof  of~\cite[1.15.1]{KP}.
 Indeed, as proved in Proposition~£3.11 above, we have a {strict covering\/} of
 $Q\-$interior $k\-$algebras
 $$k\otimes_\O {\cal B}_\delta\too k\otimes_\O S_\delta
 \leqno  £3.16.1$$
 but {\it not\/} a  {strict covering\/} $k\otimes_\O {\cal B}\too k\otimes_\O S$
 of $H\-$interior $k\-$algebras as stated in~\cite[1.15]{KP}; however, it follows from
 \cite[~2.14.4 and Lemma~9.12]{P6} that the isomorphism ${\cal B}_\delta (Q_\delta)\cong S_\delta (Q)$ induced by homomorphism~£3.16.1
 \cite[4.14]{P4} forces the {embedding\/} ${\cal B}(Q_\delta)\to
 S (Q_{\bar\delta})$ where $\bar\delta$ denotes the local point of $Q$ over $S$ determined by $\delta\,;$ hence, we still have the isomorphism
 \cite[1.15.5]{KP} which allows us to complete the argument.

\vskip 1cm \noindent{\bf\large 4. Extensions of Glauberman correspondents of blocks}

\bigskip In this section, we continue to use the notation in Paragraph 3.1, namely $\O$ is a complete discrete valuation ring with an algebraically closed residue field $k$ of characteristic $p\,;$ moreover we assume that its quotient field
$\K$ has characteristic 0 and is big enough for all finite groups that we will consider; this assumption is kept throughout the rest of this paper.

\bigskip\noindent{\bf 4.1.}\quad Let $A$ be a cyclic group of order $q$, where $q$ is a power of a prime. Assume that
$G$ is an $A$-group, that $H$ is an $A$-stable normal subgroup of~$G$ and
that $b$ is $A$-stable. Note that, in this section, $b$ is not necessarily nilpotent.  Assume that $A$ and $G$ have coprime orders; by \cite[Theorem
1.2]{P5}, $G$~acts transitively on the set of all defect groups of
$G_\alpha$ and, obviously, $A$ also acts on this set; hence, since $A$ and
$G$ have coprime orders, by \cite[Lemma 13.8 and Corollary 13.9]{I}
$A$ stabilizes some   defect group of $G_\alpha$  and $G^A$ acts
transitively on the set of them. Similarly, $A$ stabilizes some defect group
 of $N_\beta$ and $N^A$ acts transitively on the set of them. Thus, we may
assume that $A$ stabilizes $P\i N$ and actually we ssume that $A$
centralizes~$P\,;$ recall that~$Q=P\cap H$.

\bigskip\noindent{\bf 4.2.}\quad Clearly
$H^A$ is normal in $G^A$. We claim that $N^A$ is the stabilizer of
${\it w}(b)$ in $G^A$. Indeed, for any $x\in G^A$, $b^x$ is a block
of $H$ and $Q^x$ is a defect group of $b^x$; since $A$ stabilizes
$b^x$ and centralizes~$Q^x$, ${\it w}(b^x)$ makes sense. Note that
$G$ acts on ${\rm Irr}_\K(H)\,,$ that $G^A$ acts on ${\rm Irr}_\K(H^A)$
and that the Glauberman correspondence $\pi(G, A)$ is compatible
with the obvious actions of $G^A$ on ${\rm Irr}_\K(H)$ and ${\rm
Irr}_\K(H^A)$. So we have
\begin{eqnarray*}
{\rm Irr}_\K(H^A, {\it w}(b^x)) &=& \pi(H,
A)({\rm Irr}_\K(H, b^x)) \\
&=& \pi(H, A)({\rm Irr}_\K(H,b)^x )\\
&=& (\pi(H, A)({\rm Irr}_\K(H, b)))^x \\
&=& {\rm Irr}_\K(H^A,{\it w}(b))^x\,;
\end{eqnarray*}
in particular, we get ${\it w}(b^x)={\it w}(b)^x$ and therefore we have ${\it w}(b)^x={\it w}(b)$ if and only if $x$ belongs to $N^A$. We set
\begin{center}${\it w}(c)={\rm Tr}^{G^A}_{N^A}({\it w}(b))$,
${\it w}(\beta)=\{{\it w}(b)\}$ and ${\it w}(\alpha)=\{{\it
w}(c)\}$\quad.
\end{center}
Then ${\it w}(\beta)$ is a point of $N^A$ on $\O (H^A)$, ${\it w}(\alpha)$ is a point of $G^A$ on $\O (H^A)$, we have
$(N^A)_{{\it w}(\beta)}\leq (G^A)_{{\it w}(\alpha)}$ and any defect
group of $(N^A)_{{\it w}(\beta)}$ is a defect group of $(G^A)_{{\it
w}(\alpha)}$.

\bigskip\noindent{\bf 4.3.}\quad Let ${\frak B}$ and ${\it w}({\frak B})$ be the
respective sets of  $A$-stable blocks of $G$ covering $b$ and of $G^A$ covering
${\it w}(b)\,.$ Take $e\in {\frak B}\,;$ since $P$ is a defect
group of $G_\alpha$ and $c$ fulfills $ec=e$, $e$ has a defect group
contained in $P$ and therefore, since $A$ centralizes  $P\,,$ $e$ has a defect
group centralized by~$A\,;$ hence, by \cite[Proposition 1 and Theorem 1]{W}, ${\it w}(e)$ makes sense and $A$ stabilizes all the characters in ${\rm Irr}_\K(G, e)\,;$ that is to say, $A$ stabilizes all the characters of blocks in ${\frak B}$. Moreover, by \cite[Theorem 13.29]{I}, ${\it w}(e)$ belongs to ${\it w}({\frak B})$.

\bigskip\noindent{\bf Proposition 4.4.}\quad {\it The map ${\it w}:
{\frak B}\rightarrow {\it w}({\frak B}),\, e\mapsto {\it w}(e)$ is bijective and we have
\begin{center}${\rm Irr}_\K(G^A, {\it w}(c))=\pi(G,\,A)({\rm
Irr}_\K(G, c)^A)$\quad .\end{center} }

\smallskip\noindent{\it Proof.}\quad Assume that $g\in {\frak B}$ and ${\it w}(e)={\it w}(g)\,;$ then there exist $\chi\in {\rm Irr}_\K(G, e)$ and $\phi\in
{\rm Irr}_\K(G, g)$ such that $\pi(G, A)(\chi)=\pi(G, A)(\phi)$; but
this contradicts the bijectivity of the Glauberman correspondence.
Therefore the map ${\it w}$ is injective.

Take $h\in {\it w}({\frak B})\,;$ then $h$ covers ${\it w}(b)$ and so there
exist $\zeta\in {\rm Irr}_\K(G^A, h)$ and $\eta\in {\rm Irr}_\K(H^A,
{\it w}(b))$ such that $\eta$ is a constituent of ${\rm
Res}^{G^A}_{H^A}(\zeta)$. Set
$$\theta=(\pi(G, A))^{-1}(\zeta)\qq \vartheta=(\pi(H, A))^{-1}(\eta)
\quad ;$$
 by \cite[Theorem 13.29]{I}, $\vartheta$ is a constituent of
 ${\rm Res}^G_H(\theta)\,;$ let $l$ be
the block of~$G$  acting as the identity map on a $\K
G$-module affording $\theta\,;$ then $l$ covers $b$ and we have ${\it
w}(l)=h$. Finally, we have
\begin{eqnarray*}
\pi(G,\,A)({\rm Irr}_\K(G, c)^A) &=& \pi(G,\,A)(\cup_{e\in {\frak B}}{\rm Irr}_\K(G, e))  \\
   &=& \cup_{{\it w}(e)\in {\it w}({\frak B})}{\rm Irr}_\K(G^A, {\it w}(e)) \\
  &=& {\rm Irr}_\K(G^A, {\it w}(c))
\end{eqnarray*}

\bigskip\noindent{\bf Proposition 4.5.}\quad {\it $P$ is a defect group of the pointed
group $(G^A)_{{\it w}(\alpha)}$.}

\smallskip\noindent{\it Proof.}\quad It suffices  to show
that $P$ is a defect group of~$(N^A)_{{\it w}(\beta)}$ (cf.~£3.1); thus, without loss
of generality, we can assume that $G=N$. Obviously, $A$ stabilizes
$P\.    H$ and $b$ is the unique block of $P\.     H$ covering the block $b$ of~$H\,;$ since $P$ is a defect group of $G_\alpha$ and $N_\beta$, $P$ is
maximal in $N$ such that ${\rm Br}_P^{\O H}(b)\neq 0\,;$ thus $P$ is
maximal in $P\.     H$ such that ${\rm Br}_P^{\O (P\.     H)}(b)\neq 0\,;$ therefore
$P$ is a defect group of $b$ as a block of $P\.    H$. Since $A$ centralizes $P$, the
Glauberman correspondent $b'$ of $b$ as a block of $P\.     H$ makes
sense; moreover by Proposition 4.4, $b'$ covers ${\it w}(b)$. Since
${\it w}(b)$ is the unique block of $P\.     H^A$ covering the block ${\it
w}(b)$ of $H^A$, $b'={\it w}(b)$, and then, by \cite[Theorem 1]{W},
$P$ is a defect group of ${\it w}(b)$ as a block of $P\.     H^A$; in particular, ${\rm Br}_P^{\O
(H^A)}({\it w}(b))\neq 0$.

Since $P$ is a defect group of $G_\alpha$, by \cite[Theorem
5.3]{KP} the image of $P$ in the quotient group $N/H$ is a Sylow
$p$-subgroup of $N/H\,;$  but, the inclusion map $N^A\hookrightarrow N$
induces a group isomorphism $N^A/H^A\cong (N/H)^A\,;$ hence, the
image of $P$ in $N^A/H^A$ is a Sylow $p$-subgroup of $N^A/H^A\,;$ then, by
\cite[Theorem 5.3]{KP} again, $P$ is a defect group of $(N^A)_{{\it
w}(\alpha)}$.

\bigskip\noindent{\bf 4.6.}\quad We may assume that $A$ stabilizes $P_\gamma\,;$ then $A$ stabilizes $Q_\delta$ too (see \cite[Proposition 5.5]{KP}). Let $R$ be a subgroup such that $Q\leq R\leq P$ and $R_\varepsilon$  a local
pointed group on $\O H$ contained in~$P_\gamma$. Since $A$ stabilizes $P_\gamma$ and centralizes $P$,
$A$ centralizes $R$ and then, by \cite[Proposition 5.5]{KP}, it stabilizes~$R_\varepsilon$.
Since~${\rm Br}_R^{\O H}(\varepsilon)$ is a point of
$kC_H(R)$, then there is a unique block $b_\varepsilon$ of $\O C_H(R)$
such that ${\rm Br}_R^{\O H}(b_\varepsilon\varepsilon)={\rm
Br}_R^{\O H}(\varepsilon)$ and, by \cite[Lemma 2.3]{Z},
$C_Q(R)$ is a defect group of $b_\varepsilon$; in particular,
$b_\varepsilon$ is nilpotent. Obviously, $A$
centralizes $C_Q(R)$ and, since $A$
stabilizes $R_\varepsilon$ and thus it stabilizes $b_\varepsilon$, ${\it
w}(b_\varepsilon)$ makes sense; moreover, ${\it w}(b_\varepsilon)$ is
nilpotent and, since we have
$$C_{H^A}(R)=C_{C_H(R)}(A)\quad ,$$
 there
is a unique local point ${\it w}(\varepsilon)$ of $R$ on $\O (H^A)$
such that
\begin{center}${\rm Br}_R^{\O (H^A)}({\it w}(\varepsilon){\it w}(b_\varepsilon))={\rm
Br}_R^{\O (H^A)}({\it w}(\varepsilon))\quad .$\end{center}

\bigskip\noindent{\bf Proposition 4.7.}\quad {\it $P_{{\it w}(\gamma)}$ is a defect
pointed group of $(G^A)_{{\it w}(\alpha)}$ and $Q_{{\it w}(\delta)}$
is a defect pointed group of $(H^A)_{\{{\it w}(b)\}}$.}

\smallskip\noindent{\it Proof.}\quad
By \cite[Proposition 2.8]{P7}, the inclusion map $\O H\hookrightarrow \O
(P\.H)$ is actually a strict semicovering $P\.     H$-algebra homomorphism; hence, $\gamma$ determines a unique
local point $\gamma'$ of $P$ on $\O (P\.     H)$ such that $\gamma \subset
\gamma'$. Obviously, $b$ is a block of $P\.     H$. Since $\beta$ is also a
point of $P\.     H$ on $\O H$ and $P_\gamma$ is also a defect pointed
group of the pointed group $(P\.     H)_\beta$ on $\O H$, by \cite[Corollary 6.3]{KP} $P_{\gamma'}$
is a defect pointed group of the pointed
group $(P\.     H)_\beta$ on $\O (P\.     H)$.

Let $b_{\gamma'}$ be the block of $C_{P\.     H}(P)$ such that
$${\rm Br}_P^{\O (P\.     H)}(b_{\gamma'}\gamma')={\rm
Br}_P^{\O (P\.     H)}(\gamma')\quad ;$$
 then $Z(P)$ is a defect group of
$b_{\gamma'}$ and therefore ${\it w}(b_{\gamma'})$ makes sense. Obviously,
$b_{\gamma'}$ covers $b_\gamma$ and thus ${\it w}(b_{\gamma'})$
covers ${\it w}(b_{\gamma})$ (see Proposition 4.4); but, since ${\it w}(b)$
is also the Glauberman correspondent of the block $b$ of $P\.     H$ (see
the first paragraph of the proof of Proposition 4.5), by \cite[Proposition 4]{W} we have
 \begin{center}${\rm Br}_P^{\O (P\.     H^A)}({\it w}(b){\it
w}(b_{\gamma'}))={\rm Br}_P^{\O (P\.     H^A)}({\it w}(b_{\gamma'}))$
\quad ;\end{center}
this forces ${\rm Br}_P^{\O (H^A)}({\it w}(b){\it
w}(b_{\gamma}))={\rm Br}_P^{\O (H^A)}({\it w}(b_{\gamma}))$, which
implies that
$$P_{{\it w}(\gamma)}\leq (P\.     H^A)_{{\it w}(\beta)}\leq
(G^A)_{{\it w}(\alpha)}\quad ;$$
hence, by Proposition 4.5, $P_{{\it w}(\gamma)}$ is a
defect pointed group of $(G^A)_{{\it w}(\alpha)}$.

The statement that $Q_{{\it w}(\delta)}$ is a defect pointed group
of $(H^A)_{\{{\it w}(b)\}}$ is clear.

\bigskip\noindent{\bf Lemma 4.8.}\quad {\it Let $R_\varepsilon$ and $T_\eta$ be
local pointed groups on $\B$ such that $R$ is normal in $T$ and that we have
 $Q_\delta\leq R_\varepsilon\leq P_\gamma$ and  $Q_\delta\leq T_\eta\leq P_\gamma$. Then,
 we have $R_\varepsilon\leq T_\eta$ if and only if we have $${\rm Br}_T^{\O
C_H(R)}(b_\eta b_\varepsilon)={\rm Br}_T^{\O C_H(R)}(b_\gamma)\quad .$$}

\par\noindent{\it Proof.}\quad Obviously, $\B$ is a
$p$-permutation $P\.    H$-algebra (see \cite[Def. 1.1]{BP1}) by
$P\.    H$-conjugation and $(T, {\rm Br}_T^{\B}(b_{\eta}))$ and
$(R, {\rm Br}_R^{\B}(b_{\varepsilon}))$ are   $(b, P\.    H)\-$Brauer pairs
(see \cite[Def. 1.6]{BP1}). Moreover $T$ stabilizes $b_\varepsilon$, and
$\eta$ and $\varepsilon$ are the unique local points of $T$ and $R$ on $\B$
(see \cite[Proposition 5.5]{KP}) such that
\begin{center}${\rm Br}_T^{\B}(\eta){\rm Br}_T^{\B}(b_{\eta})={\rm Br}_T^{\B}(\eta)$ and ${\rm Br}_R^{\B}(\varepsilon){\rm Br}_R^{\B}(b_{\varepsilon})={\rm Br}_R^{\B}(\varepsilon)\quad .$\end{center}

Assume that $R_\varepsilon\leq T_\eta\,;$ then, there are $h\in \eta$ and
$l\in \varepsilon$ such that $hl=l=lh\,;$ thus, we have
$${\rm Br}_R^{\B}(hl)={\rm Br}_R^{\B}(l)\qq {\rm Br}_R^{\B}(h){\rm
Br}_R^{\B}(b_{\varepsilon})\neq 0\quad .$$
 Then, it follows from \cite[Def. 1.7]{BP1} that
 $$(R, {\rm Br}_R^{\B}(b_{\varepsilon}))\subset (T,{\rm Br}_T^{\B}(b_{\eta}))$$
and from \cite[Theorem 1.8]{BP1} that we have
${\rm Br}_T^{\O C_{H}(R)}(b_{\eta} b_{\varepsilon})={\rm Br}_T^{\O
C_{H}(R)}(b_{\eta})$.

Conversely, if we have
$${\rm Br}_T^{\O C_{H}(R)}(b_{\eta} b_{\varepsilon})={\rm Br}_T^{\O C_{H}(R)}(b_{\eta})$$
then, by \cite[Theorem 1.8]{BP1} we still have
 ${\rm Br}_R^{\B}(e_{\varepsilon} h)={\rm Br}_R^{\B}(h)$ for any $h\in \eta$; hence,  by the lifting theorem for idempotents, we get $R_{\varepsilon}\leq T_\eta$.

\bigskip Let $\cal R$ be a Dedekind domain of characteristic 0,
$\pi$ be a finite set of prime numbers such that $l\cal R\neq \cal
R$ for all $l\in \pi\,,$ and   $X$ and $Y$ be finite groups with $X$
acting on $Y$. We consider the group algebra ${\cal R} Y$ and set
$$Z_{\rm id}({\cal R} Y)=\oplus {\cal R} c$$ where $c$ runs over all
central primitive idempotents of ${\cal R}Y$. Obviously, $X$ acts on
$Z_{\rm id}({\cal R} Y)$ and, in the case that $X$ is a solvable $\pi$-group, Lluis
Puig exhibes a $\cal R$-algebra homomorphism ${\mathcal G l}^Y_X:
Z_{\rm id}({\cal R} X)\rightarrow Z_{\rm id}({\cal R} Y^X)$ (see
\cite[Theorem 4.6]{P7}), which unifies the usual Brauer homomorphism and
the Glauberman correspondence of characters --- called
the {\it Brauer-Glauberman correspondence\/}.

\bigskip\noindent{\bf Proposition 4.9.}\quad {\it Let $R_\varepsilon$ and $T_\eta$ be  local pointed groups on $\B$ such that $Q_\delta\leq R_\varepsilon\leq P_\gamma$ and that $Q_\delta\leq T_\eta\leq P_\gamma$. Then $R_\varepsilon\leq T_\eta$ and $R_{{\it w}(\varepsilon)}\leq T_{{\it w}(\eta)}$ are equivalent to
each other. }

\medskip\noindent{\it Proof.}\quad By induction we can assume that $R$ is normal
and maximal in $T$; in particular, the quotient $T/R$ is cyclic. In
this case, it follows from Lemma 4.8 that the inclusion $R_{{\it
w}(\varepsilon)}\leq T_{{\it w}(\eta)}$ is equivalent to
$${\rm Br}_T^{\O C_{H^A}(R)}({\it w}(b_\varepsilon){\it w}(b_\eta))
={\rm Br}_T^{\O C_{H^A}(R)}({\it w}(b_\eta))
\quad .\leqno £4.9.1$$

Let ${\mathbb Z}$ be the ring of all rational integers and $S$ be the
complement set of $p{\mathbb  Z}\cup q\mathbb  Z$ in $\mathbb  Z\,;$
then $S$ is a multiplicatively closed set in $\mathbb  Z$. We take
the localization $S^{-1}\mathbb  Z$ of $\mathbb  Z$ at $S$ and
regard it as a subring of $\K\,;$ since we assume that $\cal K$ is big enough for
all finite groups we consider, we can assume that $\cal K$ contains an $|H|$-th
primitive root $\omega$ of unity and we set
\begin{center}${\cal
R}=(S^{-1}{\mathbb Z})[\omega]\quad .$\end{center}
Then $\cal R$ is a
Dedekind domain (see \cite[Example 2 in Page 96 and Exercise 1 in
Page 99]{AM}) and given a prime $l$, we have $l\cal R\neq \cal R$ if and
only if $l=p$ or $l=q$. We consider the group algebra ${\cal
R}C_H(R)$ and the obvious action of $(T\times A)/R\cong (T/R)\times A$ on it.

Since $\cal R$ contains an $|H|$-th primitive unity root $\omega$, the blocks
$b_\varepsilon$, $b_\eta$, ${\it w}(b_\varepsilon)$ and ${\it w}(b_\eta)$ respectively belong to
$$Z_{\rm id}({\cal R}C_H(R))\;\;,\;\; Z_{\rm id}({\cal R}C_H(T))\;\;,\;\;
Z_{\rm id}({\cal R}C_{H^A}(R))\qq Z_{\rm id}({\cal R}C_{H^A}(T))
$$
(see \cite[Charpter IV, Lemma 7.2]{F}); then, by \cite[Corollary 5.9]{P7}, we have
$${\mathcal G l}^{C_{H}(R)}_{A}(b_\varepsilon)={\it w}(b_\varepsilon)\qq {\mathcal G l}^{C_{H}(T)}_{A}(b_\eta)={\it w}(b_\eta)
\quad .$$

If $R_\varepsilon\leq T_\eta$, by Lemma 4.8 we have the
equality
$${\rm Br}_{T}^{\O C_H(R)}(b_\varepsilon b_\eta)={\rm
Br}_{T}^{\O C_H(R)}(b_\eta)$$
 which is equivalent to
${\mathcal G l}^{C_H(R)}_{T/R}(b_\varepsilon)b_\eta=b_\eta$ (see
\cite[4.6.1 and the proof of Corollary 3.6]{P7}). Then by
\cite[4.6.2]{P7}, we have
\begin{eqnarray*} {\it w}(b_\eta) &=& {\mathcal G l}^{C_H(T)}_{A}(b_\eta)= {\mathcal G l}^{C_H(T)}_{A}({\mathcal G l}^{C_H(R)}_{T/R}
(b_\varepsilon)b_\eta) \\&=&
{\mathcal G l}^{C_H(R)}_{(T/R)\times A}(b_\varepsilon)
{\mathcal G l}^{C_H(T)}_{A}(b_\eta) \\
&=&{\mathcal G l}^{C_{H^A}(R)}_{T/R}
({\mathcal G l}^{C_{H}(R)}_{A}(b_\varepsilon))
{\mathcal G l}^{C_H(T)}_{A}(b_\eta)\\
&=&{\mathcal G l}^{C_{H^A}(R)}_{T/R} ({\it w}(b_\varepsilon)) {\it w}(b_\eta)\quad. \end{eqnarray*}
which is equivalent again to equality~£4.9.1 above (see \cite[4.6.1 and the proof of Corollary 3.6]{P7} and therefore it implies $R_{{\it w}(\varepsilon)}\leq T_{{\it
w}(\eta)}$. The prove that $R_{{\it w}(\varepsilon)}\leq T_{{\it
w}(\eta)}$ implies $R_\varepsilon\leq T_\eta$ is similar.

\bigskip\noindent{\bf 4.10.}\quad The assumptions and consequences above are very scattered; we collect them in this paragraph, so that readers can easily find them and we can conveniently quote them later.
Let $A$ be a cyclic group of order~$q$, where $q$ is a prime number; we
assume that $G$ is an $A$-group, that $H$ is an $A$-stable normal
subgroup of $G\,,$ that $b$ is $A$-stable, that $A$ centralizes
$P$ and stabilizes $P_\gamma\,,$ and  that $A$ and $G$ have coprime
orders. Without loss of generality, we may assume that $P\leq N$.
Then, $A$ centralizes $Q$ and stabilizes $Q_\delta\,,$ so that the Glauberman
correspondent ${\it w}(b)$ of the block $b$ makes sense; moreover, the block
${\it w}(b)$ determines two pointed group $(N^A)_{{\it w}(\beta)}$
and $(G^A)_{{\it w}(\alpha)}$ such that $(N^A)_{{\it w}(\beta)}\leq
(G^A)_{{\it w}(\alpha)}$ (see them in Paragraph 4.2)., and the local
pointed groups $P_\gamma$ and $Q_\delta$ determine respective defect pointed
groups $P_{{\it w}(\gamma)}$ and $Q_{{\it w}(\delta)}$ of
$(G^A)_{{\it w}(\alpha)}$ and $(H^A)_{{\it w}(\beta)}$ (see
Paragraph 4.6 and Proposition 4.7); actually, by Proposition 4.9, we have
$Q_{{\it w}(\delta)}\leq P_{{\it w}(\gamma)}$. Take ${\it w}(i)\in
{\it w}(\gamma)$ and ${\it w}(j)\in {\it w}(\delta)\,,$ and set
\begin{center} $(\O G^A)_{{\it w}(\gamma)}={\it w}(i)(\O G^A){\it
w}(i)$ , $(\O H^A)_{{\it w}(\gamma)}={\it w}(i)(\O H^A){\it w}(i)$\\
 and
  $(\O H^A)_{{\it w}(\delta)}={\it w}(j)(\O H^A){\it w}(j)\quad ;$\end{center}
then, $(\O G^A)_{{\it w}(\gamma)}$ is a $P$-interior  and
$(N^A/H^A)$-graded algebra; moreover, the $Q$-interior algebra $(\O
H^A)_{{\it w}(\delta)}$ with the group homomorphism
$$Q\too (\O H^A)_{{\it w}(\delta)}^*\quad , \quad u\mapsto u{\it w}(j)$$
 is a source
algebra of the block algebra $\O H^A{\it w}(b)$ (see \cite{P5}).

\vskip 1cm \noindent{\bf\large 5. A Lemma}

\bigskip From now on, we  use the notation and assumption in Paragraphs 3.1, 3.2 and 4.10; in particular, we assume that the block $b$ of $H$ is nilpotent. Obviously, $N_G(Q_\delta)$ acts on ${\rm Irr}_\K(H, b)$ and
${\rm Irr}_\K(Q)$ via the corresponding conjugation conjugation. Since $b$  is nilpotent, there is an explicit bijection between ${\rm Irr}_\K(H, b)$
and ${\rm Irr}_\K(Q)$ (see \cite[Theorem 52.8]{T}); in this section, we
will show that this bijection is compatible with the
$N_G(Q_\delta)$-actions; our main purpose is to obtain Lemma 5.6 below as a
consequence of this compatibility.

\bigskip\noindent{\bf 5.1.}\quad For any $x\in N_G(Q_\delta)$, $xjx^{-1}$ belongs to $\delta$ and thus there is some
invertible element $a_x\in \B^Q$ such that $xjx^{-1}=a_x ja_x^{-1}\,;$ let us
denote by $X$ the set of all elements $(a_x^{-1}x)j$ such that $a_x$ is invertible in $\B^Q$ and we have $xjx^{-1}=a_x
ja_x^{-1}$ when $x$
runs over $N_G(Q_\delta)$. Set
\begin{center}$E_G(Q_\delta)= N_G(Q_\delta)/QC_H(Q)\quad ;$\end{center}
then, the following equality
$$\Big((a_x^{-1}x)j\Big)\.     \Big((a_y^{-1}y)j\Big)=
\Big((a_x^{-1}xa_y^{-1}x^{-1})xy\Big)j$$
shows that $X$ is a group with respect to the multiplication and it is easily checked
that $Q\.     (\B_\delta^Q)^*$ is normal in $X$ and that the map
$$E_G(Q_\delta)\too X/Q(\B_\delta^Q)^*\leqno 5.1.1$$ sending the coset of $x\in
N_G(Q_\delta)$ in $N_G(Q_\delta)/QC_H(Q)$ to the coset of
$(a_x^{-1}x)j$ in $X/Q(\B_\delta^Q)^*$ is a group isomorphism.

\bigskip\noindent{\bf 5.2.}\quad
We denote by $Y$ the set of all such elements $a_x^{-1}x$ when $x$
runs over $N_G(Q_\delta)$ and $a_x$ over the invertible element of
$\B^Q$ such that $a_x^{-1}x$ commutes with $j$. As in 5.1, it is
easily checked that $Y$ is a group with respect to the
multiplication
$$(a_x^{-1}x)\.     (a_y^{-1}y)= (a_x^{-1}xa_y^{-1}x^{-1})xy \quad,$$
that $Y$ normalizes $Q\.     ((\O H)^Q)^*$ and that the map
$$E_G(Q_\delta)\too
\Big(Y\.     Q\.     (\B^Q)^*\Big)\Big/\Big(Q\.     (\B^Q)^*\Big)
\leqno 5.2.1$$ sending the
coset of $x\in N_G(Q_\delta)$ to the coset of $a_x^{-1}x$ in the right-hand quotient is a group isomorphism.

\bigskip\noindent{\bf 5.3.}\quad Let $I$ and $J$ be the sets of isomorphism
classes of all simple $\K\otimes_\O \B\-$ and $\K\otimes_\O \B_\delta\-$modules respectively. Cleraly, $Y$ acts on $I\,;$ but, since $Y\cap (\B^Q)^*$ acts trivially on $I$, the action of $Y$ on $I$ induces an
action of $E_G(Q_\delta)$ on $I$ through isomorphism 5.2.1;
actually, this action coincides with the action of $E_G(Q_\delta)$ on ${\rm
Irr}_\K(H, b)$ induced by the $N_G(Q_\delta)$-conjugation. Similarly,
$X$ acts on $J$ and this action of $X$ on $J$ induces an action of
$E_G(Q_\delta)$ on $J$ through isomorphism 5.1.1. But, by
\cite[Corollary 3.5]{P5}, the functor $M\mapsto j\. M$ is an equivalence
between the categories of finitely generated
$\B$- and $\B_\delta$-modules, which induces a bijection between the sets $I$ and $J$. Then, since $Y$ commutes with $j$ and the
map
$$Y\too X\quad ,\quad y\mapsto yj$$
 is a group homomorphism, it is
easily checked that this bijection is com-patible with the actions of
$E_G(Q_\delta)$ on $I$ and $J$.

\bigskip\noindent{\bf 5.4.}\quad Recall that (cf.~£3.7)
$$\B_\delta\cong T\otimes_\O \O Q
\leqno 5.4.1$$
where $T = {\rm End}_\O (W)$ for an
endo-permutation $\O Q$-module $W$ such that the determinant of the
image of any element of $Q$ in  is one; in this case, the $\O Q$-module $W$ with these properties is unique up to isomorphism. Then, for
any simple $\K\otimes_\O \B_\delta$-module $V$ there is a $\K
Q$-module $V_W$, unique up to isomorphism, such that
$$V\cong W\otimes_\O V_W$$ as $\K\otimes_\O \B_\delta$-modules; moreover the correspondence
$$V\mapsto V_W\leqno 5.4.2$$
determines a bijection between $J$ and the set of isomorphism classes of all
simple $\K Q$-modules. Now, the composition of this bijection with the bijection between isomorphism classes in 5.3 is a bijection from $I$ to the set of isomorphism classes of all simple $\K Q$-modules; translating this bijection to characters,
we obtain a bijection
$${\rm Irr}_\K(H, b)\too {\rm Irr}_\K (Q)\quad ,\quad
\chi_\lambda\mapsto \lambda
\quad ;\leqno 5.4.3$$
let us denote by $\chi\in {\rm Irr}_\K(H, b)$ the image ofthe trivial character of
$Q\,.$

\bigskip\noindent{\bf 5.5.}\quad Moreover, the $N_G(Q_\delta)$-conjugation induces an action of $E_G(Q_\delta)$ on the set of isomorphism classes of all simple $\K Q$-modules and we claim that,  for any simple $\K\otimes_\O \B_\delta$-module $V$
and any $\bar x\in E_G(Q_\delta)\,,$ we have a  $\K Q$-module isomorphism
$$^{\bar x}(V_W)\cong (^{\bar x}V)_W
\quad ;\leqno 5.5.1$$
in particular, bijection 5.4.2 is compatible with the
actions of $E_G(Q_\delta)$ on $J$ and on the set of
isomorphism classes of  simple $\K Q$-modules. Indeed, let $x$ be a lifting of
$\bar x$ in $N_G(Q_\delta)$ and denote by
$\varphi_x$ the isomorphism
$$Q\cong Q\quad ,\quad u\mapsto xux^{-1}
\quad ;$$
take a lifting $y=a_x^{-1}xj$ of $\bar x$ in $X$ through isomorphism 5.1.1;
since the conjugation by $y$ stabilizes $\B_\delta$, the map
$$f_y: \B_\delta\cong {\rm Res}_{\varphi_x}(\B_\delta)\quad ,\quad
a\mapsto yay^{-1}$$
is a $Q$-interior algebra isomorphism; then,
by \cite[Corollary 6.9]{P4}, we can modify $y$ with a suitable element
of $(\B_\delta^Q)^*$ in such a way that  $f_y$ stabilizes~$T\,;$ in this case, the restriction of $f_y$ to $T$ has to be inner and thus we have
$W\cong {\rm Res}_{f_y}(W)$ as ${\rm T}\-$modules. Moreover, since the action of $Q$ on $T$ can be uniquely lifted to a $Q$-interior algebra structure
such that the determinant of the image of any $u\in Q$ in $T$ is one, $f_y$ also stabilizes the image of $Q$ in~$T\,;$ more precisely, $f_y$ maps the image of $u\in Q$ onto
the image of $\varphi_x (u)\,.$ The claim follows.

\bigskip\noindent{\bf Lemma 5.6.}\quad {\it With the notation above,

\smallskip\noindent{\bf 5.6.1.}\quad The irreducible character $\chi$  is
$N_G(Q_\delta)$-stable and its restriction to the set $H_{p'}$ of all
$p$-regular elements of $H$ is the unique irreducible Brauer character of $H\,.$

\smallskip\noindent{\bf 5.6.2.}\quad The Glauberman correspondent
$\phi$ of $\chi$ is $N_{G^A}(Q_{{\it w}(\delta)})$-stable and its restriction
 to the set $H^A_{p'}$ of all
$p$-regular elements of $H^A$ is the unique irreducible Brauer character of $H^A\,.$}

\smallskip\noindent{\it Proof.}\quad It follows from £5.3 and £5.5 that the
bijection~£5.4.3 is compatible with the actions of
$E_G(Q_\delta)$ in ${\rm Irr}_\K(H, b)$ and ${\rm Irr}_\K (Q)\,;$
hence, $\chi$ is $E_G(Q_\delta)$-stable and thus
$N_G(Q_\delta)$-stable. Since $\phi$ is the unique irreducible
constituent of ${\rm Res}^H_{H^A}(\chi)$ occurring with a
multiplicity coprime to $q$ and $N_{G^A}(Q_{{\it w}(\delta)})$ is contained
in $N_G(Q_\delta)$, $\phi$ has to be $N_{G^A}(Q_{{\it
w}(\delta)})$-stable. By the very definition of the bijection 5.4.3,
the restriction of $\chi$ to $H_{p'}$ is the unique Brauer character
of $H$. Since the perfect isometry $R_H^b$ between $ {\cal R}_\K (H,
b)$ and ${\cal R}_\K (H^A, {\it w}(b))$ maps $\psi\in I$ onto
$\pm\pi(H, A)(\psi)$ and the blocks $b$ and ${\it w}(b)$ are
nilpotent, by \cite[Theorem 4.11]{B} the decomposition
matrices of $b$ and ${\it w}(b)$ are the same if the characters
indexing their columns correspond to  each other by the Glauberman correspondence;
hence, the restriction of~$\phi$ to $H^A_{p'}$ is the unique
Brauer character of $H^A$.

\vskip 1cm \noindent{\bf\large 6. A $k^*$-group isomorphism
$(\skew3\hat {\bar N}^{^k})^A\cong \,\widehat{\overline{\!N^A}}^{k}$}

\bigskip\noindent{\bf 6.1.}\quad
Let $xH$ be an $A$-stable coset in $\bar N$. We consider the action
of $H\rtimes A$ on $xH$ defined by the obvious action of $A$ on $xH$
and the right multiplication of $H$ on $xH\,;$ since $A$ and $G$ have
coprime orders, it follows from \cite[Lemma 13.8 and Corollary 13.9]{I}
that $xH\cap N^A$ is non-empty and that $H^A$ acts transitively on it; consequently,
we have $\bar N^A= (H\.     N^A)/H$ and the inclusion $N^A\subset N$ induces a
group isomorphism
$$\overline{ \!N^A}\cong \bar N^A= (H\.     N^A)/H
\quad .\leqno 6.1.1
$$
Note that  if $G=H \.     G^A$ then we have $\bar N^A= \bar N\,.$

\bigskip\noindent{\bf 6.2.}\quad It follows from Lemma~£5.6 that
$N = H\. N_G (Q_\delta)$ stabilizes $\chi$ and actually the central extension
$\skew3\hat{\bar N}$ of $\bar N$ by $U$ in~£3.9 above is nothing but
the so-called {\it Clifford extension\/} of $\bar N$ over $\chi\,;$ moreover, since $A$ and $U$ also have coprime orders, we can prove as above that
$\skew3\hat{\bar N}^A$ is a  central extension of $\bar N^A$ by~$U\,,$ which is
the {\it Clifford extension\/} of $\bar N^A$ over $\chi\,.$
Since the Glauberman correspondent ${\it w}(b)$ is nilpotent, we can repeat
all the above constructions for $G^A$, $H^A\,,$ ${\it w}(b)$ and $N^A\,;$ then, denoting by $U_A$ the group of $\vert H^A\vert\-$th
roots of unity, we obtain  a central extension $\,\widehat{\overline{\!N^A}}$ of
$\,\overline{ \!N^A} =\bar N^A$ by~$U_A\,,$ which is the {\it Clifford extension\/} of $\bar N^A$ over $\phi\,;$ moreover, note that $U_A$ is contained in~$U\,.$

\bigskip\noindent{\bf 6.3.}\quad At this point, it follows from  \cite[Corollary 4.16]{P8} that there is an extension group isomorphism
$$\hat N^A\cong (U\times \,\widehat{\!N^A})/\Delta_{-1} (U_A)
\leqno £6.3.1$$
where we are setting $\Delta_{-1}(U_A) = \{(\xi^{-1},\xi)\}_{\xi\in U_A}\,;$
moreover, according to \cite[Remark 4.17]{P7}, this isomorphism is defined by a
sequence of Brauer homomorphisms --- in different characteristics --- and, in particular, it is quite clear that it maps any $y\in H\i \hat N^A$ in the classes of $(1,y)$ in the right-hand member, so that  isomorphism~£6.3.1 induces a new
 extension group isomorphism
$$\skew3\hat{\bar N}^A\cong (U\times
\,\widehat{\overline{\!N^A}})/\Delta_{-1} (U_A)
\quad .$$
Consequently, denoting by $\varpi_A\,\colon U_A\to k^*$ the restriction of
$\varpi\,,$  we get a $k^*\-$group isomorphism
\begin{eqnarray*}(\skew3\hat {\bar N}^{^k})^A &=& \Big((k^*\times \skew3\hat {\bar N})/\Delta_\varrho (U)\Big)^A \cong (k^*\times \skew3\hat {\bar N}^A)/\Delta_\varrho (U) \\
&\cong &(k^*\times \,\widehat{\overline{\!N^A}})/\Delta_{\varrho_A} (U_A)
\,= \,\widehat{\overline{\!N^A}}^k
\end{eqnarray*}
as announced.

\bigskip\noindent{\bf Remark 6.4.}\quad Note that  if $G=H\.     G^A$ then
we have $\skew3\hat {\bar N}^A= \skew3\hat {\bar N}$.

\vskip 1cm \noindent{\bf\large 7. Proofs of Theorems 1.5 and 1.6}

\bigskip\noindent{\bf 7.1.}\quad
The first statement in Theorem~£1.5 follows from Propositions~£4.4 and~£4.5.
From now on, we assume that the block $b$ of $H$ is nilpotent; thus, the
Glauberman correspondent ${\it w}(b)$  is also nilpotent and
$(\O G^A){\it w}(c)$ is an extension of the nilpotent block algebra
$(\O H^A){\it w}(b)$. This section will be devoted to comparing the
extensions $\O G c$ and $\O G^A{\it w}(c)$ of the nilpotent block
algebras $\O H b$ and $\O H^A{\it w}(b)$. Applying Theorem 3.5 to the
finite groups $G^A$ and $H^A$ and the nilpotent block ${\it w}(b)$
of $H^A$, we get a finite group $L^A$ and respective injective and surjective group
homomorphisms
$$\tau^A: P\too L^A\qq \bar\pi^A: L^A\too \,\overline{\!N^A}$$
 such that $\bar\pi^A(\tau^A(u))=\bar u$ for any $u\in P$, that
 ${\rm Ker}(\bar\pi^A)=\tau^A(Q)$ and that they induce an equivalence
 of categories
$${\cal E}_{({\it w}(b),\, H^A,\, G^A)}\cong {\cal E}_{(1,\, \tau^A(Q),\, L^A)}
\quad .$$
Similarly, we et $\widehat{ L^A}= {\rm res}_{\bar \pi^A}(\,\widehat{\overline{\!N^A}}^k)$ and denote by
$\widehat{ \tau^A}\,\colon P\to \widehat{ L^A}$ the lifting of~$\tau^A\,;$
then, by Corollary 3.15,
there is a $P$-interior full matrix algebra ${\it w}(S_\gamma)$ such that we
have an isomorphism
$$(\O (G^A))_{{\it w}(\gamma)}\cong {\it w}(S_\gamma)\otimes_{\O }
\O_*\widehat{L^A}^\circ\quad \leqno{7.1.1}$$ of both $P$-interior
and $N^A/H^A$-graded algebras.

\bigskip\noindent{\bf Lemma 7.2.}\quad {\it Assume that $G=H\.     G^A\,.$ Then we have $N=H\.     N^A\,,$  the inclusion $N^A\subset N$ induces a group isomorphism $\,\overline{\!N^A}\cong \bar N$ and there is a group isomorphism
$$\sigma: L^A\cong L$$ such that
$\sigma\circ\tau^A=\tau$ and $\bar\pi\circ\sigma=\bar\pi^A$. }

\medskip\noindent{\it Proof.}\quad
For any subgroups $R$ and $T$ of $P$ containing $Q$, let us denote by
$${\cal E}_{(b,\, H,\, G)}(R, T)\qq {\cal E}_{({\it w}(b),\, H^A,\, G^A)}(R, T)$$
the respective sets of ${\cal E}_{(b,\, H,\, G)}\-$ and ${\cal E}_{({\it
w}(b),\, H^A,\, G^A)}\-$morphisms from $T$ to $R\,;$ since $A$ acts trivially
in ${\cal E}_{(b,\, H,\, G)}(R, T)$, by \cite[Lemma 13.8 and Corollary
13.9]{I} each morphism in ${\cal E}_{(b,\, H,\, G)}(R, T)$ is
induced by some element in $N^A\,;$ moreover, if $T_\nu$ and $R_\varepsilon$ are local pointed groups contained in $P_\gamma\,,$ it follows from Proposition 4.9 that we have $T_\nu\leq (R_\varepsilon)^x$  for some $x\in N^A$ if and only if
we have $T_{{\it w}(\nu)}\leq (R_{{\it w}(\varepsilon)})^x$.
Therefore, we get
$${\cal E}_{(b,\, H,\, G)}(T, R)={\cal E}_{({\it w}(b),\, H^A,\, G^A)}(T, R)
\quad . $$

At this point, it is easy to check that $L$, $\tau$ and $\bar\pi$
fulfill the conditions in Theorem 3.5 with respect to $G^A$, $H^A$ and the nilpotent
block ${\it w}(b)$. Then this lemma follows from the uniqueness part in
Theorem 3.5.

\bigskip\noindent{\bf Lemma 7.3.}\quad {\it Assume that $G=H\.     G^A\,.$ Then there is a $k^*$-group isomorphism   $\hat \sigma: \widehat{ L^A}\cong \hat L$ lifting $\sigma$ and fulfilling $\hat\sigma\circ \widehat{\tau^A }= \hat\tau\,.$
In particular, we have
$${\rm Irr}_{\K}(G, c)={\rm Irr}_{\K}(G, c)^A\quad .$$}

\par\noindent{\it Proof.}\quad The first statement is an easy consequence of £6.3
and Lemma 7.2; then, the last equality follows from Corollary~£3.15.

\bigskip\noindent{\bf 7.4.} {\it Proof of Theorem 1.6.}\quad
Firstly we consider the case where the block $b$ of $H$ is not
stabilized by~$G\,;$ then we have an isomorphism
$${\rm Ind}^G_{N}(\O N b)\cong \O G b$$
of $\O G$-interior algebras
mapping $1\otimes a\otimes 1$ onto $a$ for any $a\in \O N b$ and an
isomorphism
$${\rm Ind}^{G^A}_{N^A}(\O (N^A ){\it w}(b))\cong \O (G^A){\it w}(b)$$
of $\O (G^A)$-interior algebras mapping $1\otimes
a\otimes 1$ onto $a$ for any $a\in \O (N^A) {\it w}(b)$. Suppose that
an $\O(N^A\times N)$-module $M$ induces a Morita equivalence from
$\O (N^A) {\it w}(b)$ to $\O N b$. Then it is easy to see that the
$\O(G^A\times G)$-module ${\rm Ind}^{G^A\times G}_{N^A\times N}
(M)$ induces a Morita equivalence from $\O Gc$ to $\O (G^A) {\it w}(c)\,.$ So, we can assume that $G= N$ and then we have $G^A=N^A\,.$

\smallskip
By Corollary 3.15, there exists an isomorphism of both
$(N/H)$-graded  and $P$-interior algebras
$$(\O G)_\gamma\cong S_\gamma\otimes_{\O } \O_*\hat L^{^\circ}\quad ;
\leqno{7.4.1}$$
denote by $V_\gamma$ an $\O P$-module such that ${\rm End}_\O(V_\gamma)\cong S_\gamma\,;$ choosing $i\in \gamma$ and assuming that $(\O G)_\gamma
= i(\O G)i\,,$ we know that the $\O Gb\otimes_\O (\O G)_\gamma^\circ\-$module
$(\O G)i$ determines a Morita equivalence from $\O Gb$ to $(\O G)_\gamma\,,$
whereas the $(\O G)_\gamma\otimes_\O \O_*\hat L\-$module
$V_\gamma\otimes_\O \O_*\hat L^{^\circ}$ determines a Morita equivalence
from  $(\O G)_\gamma$ to~$\O_*\hat L^{^\circ}\,,$ so that the
$\O Gb\otimes_\O \O_*\hat L\-$ module
$$(\O G)i\otimes_{(\O G)_\gamma} (V_\gamma\otimes_\O \O_*\hat L^{^\circ})
\cong (\O G)i\otimes_{S_\gamma} V_\gamma$$
 determines a Morita equivalence from  $\O Gb$ to~$\O_*\hat L^{^\circ}\,.$

 \smallskip
 Similarly, choosing $j\in \delta$ such that $ji = j = ij\,,$ assuming that
 $j(\O H)j = (\O H)_\delta$ and setting $j\. V_\gamma = V_\delta\,,$
 so that $S_\delta = {\rm End}_\O (V_\delta)\,,$  the
$\O Hb\otimes_\O \O Q\-$ module
$$(\O H)j\otimes_{(\O H)_\delta} (V_\delta\otimes_\O \O Q)
\cong (\O H)j\otimes_{S_\delta} V_\delta$$
 determines a Morita equivalence from  $\O Hb$ to~$\O Q\,.$

 \smallskip
 Analogously, with evident notation, the $\O (G^A) w(b)\otimes_\O
 \O_*\,\widehat{\!L^A}\-$ module
$$\O (G^A) w(i) \otimes_{w(S_\gamma)} w(V_\gamma)$$
 determines a Morita equivalence from  $\O (G^A)(b)$ to~$\O_*\,\widehat{L^A}^\circ\,,$ whereas  the $\O (H^A) w(b)\otimes_\O \O Q\-$module
$$\O (H^A)w(j)\otimes_{w(S_\delta)} w(V_\delta)$$
 determines a Morita equivalence from  $\O (H^A) w(b)$ to~$\O Q\,.$

 \smallskip
 Consequently, identifying $\,\widehat{L^A}$ with $\hat L$ through the isomorphism
 $\hat\sigma$ (cf. Lemma~£7.3), the $\O (G\times G^A)\-$module
 $$D= ((\O G)i\otimes_{S_\gamma} V_\gamma)\otimes_{\O_*\hat L}
  (w(V_\gamma)^\circ \otimes_{w(S_\gamma)} w(i) \O (G^A))$$
 determines a Morita equivalence from  $\O Gb$ to~$\O (G^A) w(b)\,,$
 whereas the $\O (H\times H^A)\-$module
 $$M = ((\O H)j\otimes_{S_\delta} V_\delta)\otimes_{\O Q }
 (w(V_\delta)^\circ \otimes_{w(S_\delta)} w(j)\O (H^A))$$
 determines a Morita equivalence from  $\O Hb$ to~$\O (H^A) w(b)\,.$

 \smallskip
 Moreover, since we have the obvious inclusions
 $$(\O H)j\subset (\O G)i\quad ,\quad S_\delta \subset S_\gamma\qq
 V_\delta \subset V_\gamma
 \quad ,$$
it is easily checked that we have
$$(\O H)j\otimes_{S_\delta} V_\delta\cong (\O H)i\otimes_{S_\gamma} V_\gamma\subset (\O G)i\otimes_{S_\gamma} V_\gamma
\quad ;\leqno £7.4.2$$
in particular, we have an evident section
$$(\O G)i\otimes_{S_\gamma} V_\gamma\too (\O H)j\otimes_{S_\delta} V_\delta$$
which is actually an $\O Hb\otimes_\O \O Q\-$module homomorphism. Similarly,
we have a split $\O (H^A) w(b)\otimes_\O \O Q\-$module monomorphism
$$\O (H^A)w(j)\otimes_{w(S_\delta)} w(V_\delta)\too
\O (G^A) w(i) \otimes_{w(S_\gamma)} w(V_\gamma)
\quad .\leqno £7.4.3$$

 \smallskip
 In conclusion, the $\O Hb\otimes_\O \O Q\-$ and
 $\O (H^A) w(b)\otimes_\O \O Q\-$module homomorphisms~£7.4.2 and~£7.4.3, together with the inclusion $\O Q\subset \O \hat L\,,$ determine
  an  $\O (H\times H^A)\-$module homomorphism
 $$M\too {\rm Res}_{H\times H^A}^{G\times G^A} (D)
 \leqno £7.4.4$$
 which actually admits a section too. Now, denoting by $K$ the inverse image in
 $G\times G^A$  of the ``diagonal'' subgroup of $(G/H)\times (G^A/H^A)\,,$
 we claim that the product by $K$ stabilizes the image of $M$ in $D\,,$ so that
 $M$ can be extended to an $\O K\-$module.

 \smallskip
 Actually, we have
 $$K = (H\times H^A)\.\Delta (N_{G^A}(Q_\delta))
 \quad ,$$
 so that it suffices to prove that the image of $M$ is stable by multiplication
 by $\Delta (N_{G^A}(Q_\delta))\,.$ Given $x\in N_{G^A}(Q_\delta)$, there are
some invertible elements $a_x\in (\O H)^Q$ and  $b_x\in (\O (H^A))^Q$
such that
$$xjx^{-1} = a_xj a_x^{-1}\qq x w(j)x^{-1} = b_x w(j )b_x^{-1}$$
and therefore $a_x^{-1}x$ and $b_x^{-1}x$ respectively centralize $j$ and
$w(j)\,,$ so that $a_x^{-1}xj$ and $b_x^{-1}x w(j)$ respectively belong to
$(\O G)_\delta$ and to $(\O G^A)_{w(\delta)}\,;$ but, according to
isomorphisms~£7.4.1 and~£7.1.1, we have $G/H\-$ and $G^A/H^A\-$gra-ded
isomorphisms
$$(\O G)_\delta\cong S_\delta\otimes_\O \O_*\hat L^\circ\qq
(\O G^A)_{w(\delta)}\cong w(S_\delta)\otimes_\O \O_*\hat L^\circ
$$
where we are setting $w(S_\delta)  = w(j)w(S_\gamma)w(j)\,.$

\smallskip
Hence, identifying  with each other both members of these isomorphisms and modifying if necessary our choice of $a_x\,,$ for some $s_x\in S_\delta\,,$
$t_x\in w(S_\delta)$ and $\hat y_x\in \hat L^\circ\,,$
 we get
$$a_x^{-1}x j = s_x\otimes \hat y_x\qq b_x^{-1}x w(j) = t_x\otimes \hat y_x
\quad .$$
Thus, setting $w(V_\delta) = w(j)w(V_\gamma)\,,$  for any
$a\in (\O H)j\,,$ any $b\in (\O H^A)w(j)\,,$  any $v\in V_\delta$
and any $w\in w(V_\delta)\,,$ in $D$ we have
\begin{eqnarray*}(x,x)\!\!\!\!&\.&\!\!\!\!(a\otimes v)\otimes (w\otimes b)
= (x a\otimes v)\otimes (w\otimes bx^{-1})\\
&=& (x ax^{-1}a_x (a_x^{-1}xj)\otimes v)\otimes
(w\otimes (w(j)x^{-1}b_x)b_x^{-1}xbx^{-1})\\
&=&(x ax^{-1}a_x \otimes s_x\.v)\. \hat y_x\otimes
\hat y_x^{-1}\.(w\.t_x^{-1}\otimes b_x^{-1}xbx^{-1})\\
&=&(x ax^{-1}a_x \otimes s_x\.v) \otimes
(w\.t_x^{-1}\otimes b_x^{-1}xbx^{-1})\quad ;
\end{eqnarray*}
since $x ax^{-1}a_x $ and $b_x^{-1}xbx^{-1}$ respectively belong to
$(\O H)j$ and $w(j)(\O H^A)\,,$ this proves our claim.

\smallskip
Finally, since homomorphism~£7.4.4 actually becomes an $\O K\-$module
homomorphism, it induces an $\O (G\times G^A)\-$module homomorphism
$${\rm Ind}_K^{G\times G^A}(M)\too D$$
which is actually an isomorphism as it is easily checked. We are done.

\bigskip The following theorem is due to Harris and Linckelmann (see \cite{H}).

\bigskip\noindent{\bf Theorem 7.5.}\quad {\it Let $G$ be an $A$-group and assume that $G$ is a finite $p$-solvable group and $A$ is a solvable group of order prime to $|G|$. Let $b $ be an $A$-stable block of $G$ over $\O$ with a defect group $P$ centralized by $A$ and denote by ${\it w}(b)$ the Glauberman correspondent of the block $b$. Then the block algebras $\O Gb$ and $\O (G^A){\it w}(b)$ are basically Morita equivalent. }

\medskip\noindent{\it Proof.}\quad By \cite[Theorem 5.1]{H}, we can assume that $b$ is a $G\rtimes A$-stable block of ${\rm O}_{p'}(G)$, where ${\rm O}_{p'}(G)$ is the maximal normal $p'$-subgroup of $G$. Clearly $b$ as a block of ${\rm O}_{p'}(G)$ is nilpotent and thus $\O Gb$ is an extension of the nilpotent block algebra $\O {\rm O}_{p'}(G) b$.
By \cite[Theorem 5.1]{H} again, ${\it w}(b)$ is a $G^A$-stable block of ${\rm O}_{p'}(G^A)$ and thus is nilpotent; thus $\O (G^A) {\it w}(b)$ is an extension of the nilpotent block algebra $\O {\rm O}_{p'}(G^A) {\it w}(b)$. By \cite[Theorem 4.1]{H}, ${\it w}(b)$ is also the Glauberman correspondent of $b$ as a block of ${\rm O}_{p'}(G)$. Then, by Theorem 1.6, the block algebras $\O Gb$ and $\O (G^A){\it w}(b)$ are basically Morita equivalent.

\bigskip The following theorem is due to Koshitani and Michler (see \cite{KG}).

\bigskip\noindent{\bf Theorem 7.6.}\quad {\it Let $G$ be an $A$-group and assume that $A$ is a solvable group of order prime to $|G|$. Let $b $ be an $A$-stable block of $G$ over $\O$ with a defect group $P$ centralized by $A$ and denote by ${\it w}(b)$ the Glauberman correspondent of the block $b$. Assume that $P$ is normal in $G$. Then, the block algebras $\O Gb$ and $\O (G^A{)\it w}(b)$ have isomorphic source algebras. }

\medskip\noindent{\it Proof.}\quad Since $P$ is normal in $G$, by \cite[2.9]{AB} there is a block $b_P$ of $C_G(P)$ such that $b={\rm Tr}^G_{G_{b_P}}(b_P)$, where $G_{b_P}$ is the stabilizer of $b_P$ in $G$. Since $A$ and $G$ have coprime orders, by \cite[Lemma 13.8 and Corollary 13.9]{I}, $b_P$ can be chosen such that $A$ stabilizes $b_P$. Since $P$ is the unique defect group of $b$, $P$ has to be contained in $G_{b_P}\,;$ then by \cite[Proposition 5.3]{KP}, the intersection $Z(P)=P\cap C_G(P)$ is the defect group of $b_P$ and, in particular, $b_P$ is nilpotent. Thus the block $\O G b$ is an extension of the nilpotent block algebra~$\O(P\. C_G(P))b_P$ and, in particular, we have $\bar N \cong E_G (P_\gamma)\,.$

The Glauberman correspondent of $b_P$ makes sense and by \cite[Proposition 4]{W},
we have
 $${\it w}(b)={\rm Tr}^{G^A}_{(G^A)_{{\it w}(b_P)}}({\it w}(b_P))
 \quad .$$
 Since ${\it w}(b_P)$ has  defect group $Z(P)$, it is also nilpotent and thus $\O (G^A){\it w}(b)$ is an extension of the nilpotent block algebra
 $\O(P\. C_{G^A}(P)){\it w}(b_P)\,;$  once again, we have
 $\bar N ^A\cong E_{G^A} (P_{w(\gamma)})\,.$

 \smallskip
 On the other hand, since $P$ is normal in $G\,,$ it follows from \cite[Proposition 14.6]{P6} that
 $$(\O G)_\gamma\cong \O_*(P \rtimes
 \hat  E_G (P_\gamma))\qq  (\O G^A)_{w(\gamma)}\cong \O_*(P \rtimes \hat  E_{G^A} (P_{w(\gamma)}))
 \quad ;$$
 but, it follows from~£6.3 that we have a $k^*\-$group isomorphism
 $$\hat  E_G (P_\gamma)\cong \hat  E_{G^A} (P_{w(\gamma)})
 \quad .$$
 We are done.

Lluis Puig

CNRS, Institut de Math\'ematiques de Jussieu

6 Av Bizet, 94340 Joinville-le-Pont, France

\smallskip puig@math.jussieu.fr

\bigskip Yuanyang Zhou

Department of Mathematics and Statistics

Central China Normal University

Wuhan, 430079

P.R. China

\smallskip zhouyy74@163.com

\end{document}